\input amssym.def
\input amssym.tex

\input xy
\xyoption{all}

\documentclass[12pt]{article}
\usepackage[all]{xy}
\usepackage[frenchb]{babel}
\font\tenmsb=msbm10
\font\sevenmsb=msbm7
\font\fivemsb=msbm5

\newfam\msbfam
\textfont\msbfam=\tenmsb
\scriptfont\msbfam=\sevenmsb
\scriptscriptfont\msbfam=\fivemsb
\def\Bbb#1{{\fam\msbfam #1}}

\makeatletter
\ifnum\@ptsize=0  \fi \ifnum\@ptsize=2
 \fi 

\newcommand\sK{{\cal K}}

\newcommand\sO{{\cal O}}

\newcommand\bR{{\Bbb R}}
\newcommand\bZ{{\Bbb Z}}
\newcommand\bC{{\Bbb C}}
\newcommand\bQ{{\Bbb Q}}

\newtheorem{theorem}{Th\'eor\`eme}[section]
\newtheorem{lemma}[theorem]{Lemme}
\newtheorem{corollary}[theorem]{Corollaire}
\newtheorem{proposition}[theorem]{Proposition}
\newtheorem{question}[theorem]{Question}
\newtheorem{re}[theorem]{Remarque}

\newtheorem{definition}[theorem]{D\'efinition}
\newtheorem{conjecture}[theorem]{Conjecture}

\newtheorem{example}[theorem]{Example}

\newenvironment{remark}{\begin{re}\em}{\end{re}}

\begin{document}
\title {QUOTIENTS R\'ESOLUBLES OU NILPOTENTS DES GROUPES DE K\" AHLER ORBIFOLDES} \author{Fr\'ed\'eric Campana }

\maketitle

\tableofcontents

\section{Introduction}

\

$\bullet$ L'objectif du pr\'esent texte est, en particulier, d'\'etendre aux {\it orbifoldes g\'eom\'etriques  finies} K\"ahl\'eriennes $(X\vert \Delta)$ les r\'esultats connus concernant les ensembles de Green-Lazarsfeld et les quotients r\'esolubles des groupes fondamentaux, ainsi que les rev\^etements universels des vari\'et\'es  K\" ahl\' eriennes compactes. 

Cette extension est motiv\'ee par la classification bim\'eromorphe des vari\'et\'es K\" ahl\'eriennes compactes et l'\'etude de leurs rev\^etements universels, qui reposent de mani\`ere essentielle sur la consid\'eration de ces orbifoldes g\'eom\'etriques $(X\vert \Delta)$, qui sont des rev\^etements ramifi\'es {\it virtuels} des vari\'et\'es K\" ahl\'eriennes compactes, et \'eliminent virtuellement les fibres multiples de fibrations  (voir  les \S\ref{ros} et \S\ref{redalge} ci-dessous pour une br\`eve description, qui fournit certaines des motivations du pr\'esent texte, et [Ca 04] et [Ca 07] pour les d\'etails). 

Bien que les objets consid\'er\'es ici soient essentiellement les m\^emes\footnote{Deligne-Mostow emploient en fait, dans [D-M 93, \S 14, pp. 137-141], le terme {\it orbifolde} pour d\'esigner essentiellement les objets consid\'er\'es ici. Leur int\'er\^et y est \'egalement centr\'e sur leur groupe fondamental, d\'efini comme ci-dessous. Je remercie F. Catanese qui m'a indiqu\'e cette r\'ef\'erence.}  que les classiques champs de Deligne-Mumford lisses, le point de vue adopt\' e ici  am\`ene naturellement \`a poser les questions consid\'er\'ees ici (ainsi que de nombreuses autres), non abord\'ees jusqu'ici dans la th\'eorie des champs.

\

Les deux questions\footnote{Voir aussi \S \ref{qq}}que nous souhaitons aborder  dans un cas tr\`es particulier sont ([Ca 07, question 12.12]):

\

1. Pour toute $(X\vert\Delta)$ lisse, avec $X$ K\" ahler compacte, $\pi_1(X\vert\Delta)$ est-il commensurable \`a $\pi_1(Y)$, pour $Y$ compacte K\" ahler (ou projective) lisse ad\'equate, de m\^eme dimension que $X$?

2. Le rev\^etement universel\footnote{D'un mod\`ele bim\'eromorphe ad\'equat, ``minimal" si $(X\vert \Delta)$ n'est pas unir\'egl\'ee.} $\widetilde{(X\vert\Delta)}$ est-il isomorphe au rev\^etement universel $\tilde{Y}$ d'une  $Y$ compacte K\" ahler \`a singularit\'es quotient, et de m\^eme dimension que $X$? Peut-on choisir $Y$ projective\footnote{Les exemples de C. Voisin de vari\'et\'es K\" ahl\'eriennes non homotopiquement projectives n'infirment en rien ces coincidences potentielles, puisque soit non-minimaux, soit unir\'egl\'es. Ces exemples montrent plutot que certaines transformations bim\'eromorphes ne sont pas r\'ealisables sur des vari\'et\'es projectives.} si le rev\^etement universel de $(X\vert \Delta)$ n'est pas recouvert par des sous-vari\'et\'es complexes compactes de dimension strictement positive?

\

 Puisque l'on a une suite exacte: $K\to \pi_1(X\vert\Delta)\to \pi_1(X)\to 1$ dans laquelle $K$ est engendr\'e par des \'el\'ements de torsion en nombre fini modulo conjugaison (mais $K$ n'est ni de torsion, ni de type fini, en g\'en\'eral, voir \S. 2.2), les questions pr\'ec\'edentes admettent une r\'eponse affirmative lorsque $\pi_1(X\vert\Delta)$ est r\'esiduellement fini ou admet un sous-groupe d'indice fini sans torsion, et donc, en particulier, lorsque ce groupe est lin\'eaire, et en particulier, virtuellement nilpotent.

Par contre, les groupes r\'esolubles de pr\'esentation finie ne sont pas n\'ecessairement r\'esiduellement finis (un contre-exemple r\'esoluble au probl\`eme du mot, et donc \`a la finitude r\'esiduelle, est construit dans [K 81]).Ces groupes forment donc la classe de groupes la plus simple sur laquelle tester les questions pr\'ec\'edentes. Les r\'esultats du pr\'esent texte fourniront, en particulier, une r\'eponse affirmative \`a ces deux questions lorsque $\pi_1(X\vert \Delta)$ est r\'esoluble. 

\

  Plus pr\'ecis\'ement, nous montrons ici, entre autres choses, que les r\'esultats connus concernant les ensembles de Green-Lazarsfeld et les quotients r\'esolubles des groupes de K\" ahler s'\'etendent sans changement aux groupes fondamentaux des orbifoldes g\'eom\'etriques K\" ahl\'eriennes\footnote{Et devraient aussi s'\'etendre, semble-t-il, aux champs de Deligne-Mumford lisses dont l'espace des modules est domin\'e par une vari\'et\'e K\" ahl\'erienne compacte. Je remercie C. Simpson pour avoir soulev\'e cette question.}. 

En particulier, pour le groupe fondamental d'une telle orbifolde:

1. Son ensemble de Green-Lazarsfeld est r\'eunion d'un nombre fini de translat\'es par des \'el\'ements de torsion de sous-tores du groupe des caract\`eres de ce groupe (\S3).

2. Ses quotients r\'esolubles sont tous virtuellement nilpotents si (et seulement si) l'orbifolde consid\'er\'ee n'a pas de morphisme sur une courbe orbifolde hyperbolique (apr\`es \' eventuel rev\^etement orbifolde-\'etale fini) (\S4).

Les quotients virtuellement nilpotents des groupes fondamentaux des orbifoldes compactes K\"ahler sans morphisme sur une courbe orbifolde hyperbolique coincident donc avec ces m\^emes quotients virtuellement r\'esolubles. Les quotients nilpotents sont classifi\'es par la suite centrale descendante correspondante, qui est d\'etermin\'ee par ses deux premiers termes, par `formalit\'e' au sens de D. Sullivan. Ces quotients sont n\'eammoins extr\^emement mal connus, m\^eme lorsque $\pi_1(X)$ est suppos\'e nilpotent sans torsion (des exemples existent de classe de nilpotence $2$, l'existence est actuellement inconnue pour des classes de nilpotence $3$ ou plus).

3. Dans une derni\`ere partie (\S 7) , nous montrons, g\'en\'eralisant les r\'esultats de [Ca 01], que si $f:X\to Y$ est une application holomorphe surjective connexe, avec $X$ compacte et K\" ahler, les classes de nilpotence $\nu_s$ des s-i\`emes quotients des s\'eries centrales (quotient\'es par leur torsion) de $\pi_1(X), \pi_1(Y\vert \Delta_f)$, et $\pi_1(X_y)$ satisfont, comme dans le cas d'un produit, une in\'egalit\'e non-archim\'edienne de la forme: $\nu_s(\pi_1(X))=max\{\nu_s(\pi_1(X_y)),\nu_s(\pi_1(Y\vert \Delta_f))\}$. Ce fait, qui repose sur des r\'esultats profonds de th\'eorie de Hodge, d\^us \`a P. Deligne et R. Hain, est en complet contraste avec le cas (non-K\" ahler) des quotients d'un groupe de Lie complexes r\'esolubles par un r\'eseau cocompact.



4. Nous montrons aussi (aux \S5 et \S6) que les ensembles de Green-Lazarsfeld et les quotients r\'esolubles du groupe fondamental coincident essentiellement pour $(X\vert \Delta)$, pour la factorisation de Stein de son image d'Albanese, et (\`a un facteur ab\'elien pr\`es) pour son ``coeur". En particulier (voir \S 6), les quotients r\'esolubles de $\pi_1(X\vert \Delta)$ sont virtuellement ab\'eliens si $(X\vert \Delta)$ est {\it sp\'eciale} (voir [Ca 07] et le \S\ref{orbspec} pour les d\'efinitions du ``coeur" et de ``sp\'eciale"). Conjecturalement (\ref{cj}), $\pi_1(X\vert \Delta)$ est lui-m\^eme virtuellement ab\'elien si $(X\vert \Delta)$ est ``sp\'eciale".

On r\'eduit ainsi l'\'etude de ces invariants, pour $(X\vert \Delta)$ arbitraire, au cas crucial o\`u $(X\vert \Delta)$ est projective de type g\'en\'eral, et o\`u son morphisme d'Albanese est g\'en\'eriquement fini sur son image. 

\

 Il serait int\'eressant de savoir si ces r\'esultats subsistent pour l'image d'Albanese (et non seulement pour sa factorisation de Stein. Voir la question \ref{qqralb}), et aussi dans quelle mesure les r\'esultats pr\'ec\'edents peuvent \^etre \'etendus aux groupes fondamentaux des orbifoldes logarithmiques (ie: des vari\'et\'es ``quasi-projectives" ou ``quasi-K\" ahler") en utilisant la th\'eorie de Hodge mixte. Certains r\'esultats sont d\'ej\`a connus dans des situations particuli\`eres de ce cas ([A 98], [AN 93], [M 78]). Voir aussi le \S \ref{qq} pour d'autres questions.

\

{\bf Remerciements:} Je voudrais remercier T. Delzant pour m'avoir communiqu\'e ses textes [De 06] et [De 07], fourni une photocopie de [B-S], et signal\'e l'oubli d'une condition dans une version initiale du pr\'esent texte.

\section{Rappels sur les orbifoldes g\'eom\'etriques et leurs rev\^etements.}

\

On rappelle ici tr\`es bri\`evement certaines des notions relatives aux {\it orbifoldes g\'eom\'etriques} introduites dans [Ca 01] et [Ca 07]. On se limite \`a celles qui sont strictement indispensables pour le pr\'esent texte. Ces derni\`eres notions ont, en fait, \'et\'e d\'ej\`a d\'efinies tout simplement sous le terme ``d'orbifoldes" dans [D-M 93, \S 14, pp. 135-141], texte que m'a signal\'e F. Catanese.

\subsection{Orbifoldes g\'eom\'etriques}\label{orbgeo} Si $X$ est un espace analytique complexe normal et connexe, un {\it diviseur orbifolde} $\Delta$ sur $X$ est un $\bQ$-diviseur effectif de la forme $\Delta:=\sum_J (1-1/m_j).D_j$, dans lequel les $m_j\geq 2$ sont, soit des entiers\footnote{Il s'agit donc ici seulement de diviseurs orbifoldes {\it entiers} dans la terminologie de [Ca 07].}, soit $+\infty$, et o\`u les $D_j$ sont des diviseurs de Weil irr\'eductibles distincts de $X$, localement finis. La r\'eunion des $D_j$ est appel\'ee le {\it support} de $\Delta$, not\'e $Supp(\Delta)$. Le couple $(X\vert \Delta)$ est une {\it orbifolde g\'eom\'etrique}, dite {\it logarithmique} (resp. {\it finie}) si $m_j=+\infty, \forall j$ (resp. si $m_j<+\infty, \forall j)$. Les orbifoldes g\'eom\'etriques interpolent donc lorsque $X$ est compacte, entre les cas {\it propre} (o\`u $\Delta=0$), et logarithmique. On dira enfin que $(X\vert \Delta)$ est {\it lisse} si $X$ est lisse et si $Supp(\Delta)$ est un diviseur \`a croisement normaux. On dira enfin que $(X\vert \Delta)$ est {\it K\"ahler} si $X$ est lisse, K\"ahl\'erienne, compacte et connexe.

Pour tout diviseur de Weil irr\'eductible de $X$, on notera $m_{\Delta}(D)$ sa {\it multiplicit\'e} dans $\Delta$, \'egale \`a $m_j$ si $D=D_j$ pour l'un (unique) des $j\in J$, et \`a $1$ sinon.

On dira que l'orbifolde $\Delta$ sur $X$ {\it divise} l'orbifolde $\Delta_1$ sur $X$ si les multiplicit\'es de $\Delta$ divisent toutes celles de $\Delta_1$ (ie: si pour tout $D$, $m_{\Delta}(D)$ divise $m_{\Delta_1}(D)$).

\

Un cas particulier classique est celui des orbifoldes de courbes, dans lequel $X=C$ est une courbe projective lisse et connexe, et $\Delta$ est alors la donn\'ee d'un nombre fini de points $p_j$ affect\'es de multiplicit\'es enti\`eres $m_j$. Le fibr\'e canonique est alors le $\bQ$-diviseur $K_{C\vert \Delta}:=K_C+\Delta$, dont le degr\'e est: $deg(K_{C\vert \Delta}):=2(g-1)+\sum_j(1-1/m_j)$. On dit que $C\vert \Delta)$ est {\it hyperbolique} si $deg(K_C+\Delta)>0$.

\

Un morphisme orbifolde\footnote{ Il s'agit des morphismes orbifoldes {\it divisibles} dans la terminologie de [Ca 07] qui consid\`ere des morphismes et orbifoldes g\'eom\'etriques plus g\'en\'eraux. On ne consid\`erera ici que les morphismes divisibles, auxquels s'appliquent toutes les consid\'erations de [Ca 04] et [Ca 07].} $f:(X\vert \Delta)\to (C\vert \Delta_C)$ d'une orbifolde g\'eom\'etrique sur une courbe orbifolde est une application holomorphe surjective et \`a fibres connexes de l'espace analytique normal connexe et compact $X$ sur la courbe projective lisse $C$ telle que pour tout $b\in Supp(\Delta_C)$, si $f^*(b)=\sum_k t_k.E_k$, on a, pour tout $k$: $m_{\Delta_C}(b)$ divise $t_k.m_{\Delta}(E_k)$.

Lorsque $dim(C)$ est arbitraire, on impose que $f(X)\subsetneq \Delta_C$, et que la condition de divisibilit\'e pr\'ec\'edente soit satisfaite pour tout diviseur irr\'eductible $b$ de $C$.

\

Si $f:X\to C$ est un morphisme surjectif, avec $X,C$ comme pr\'ec\'edemment, et si $(X\vert \Delta)$ est une structure orbifolde sur $X$, on d\'efinit la {\it base orbifolde} de $(f,\Delta)$ sur $C$, not\'ee $\Delta_{f,\Delta}$ (ou simplement $\Delta_f$ s'il n'y a pas d'ambiguit\'e sur $\Delta$), comme \'etant la plus petite des structures orbifoldes $\Delta_C$ sur $C$ telle que $f:(X\vert \Delta)\to (C\vert \Delta_C)$ soit un morphisme orbifolde. Explicitement: $m_{\Delta_f}(b)=pgcd\{t_k.m_{\Delta}(E_k)\}$ pour tout $b\in C$, avec les notations pr\'ec\'edentes.

Lorsque $dim(C)>1$, la d\'efinition est la m\^eme, mais en omettant les diviseurs $f$-exceptionnels de $X$.

\begin{remark}\label{d'} Remarquons que si $d:Y\to X$ et $f:X\to C$ sont des fibrations, et si $F:=f\circ d$, alors: $\Delta_{F}$ divise $\Delta_f$, et $\Delta_f$ divise $\Delta_{f,\Delta}$ pour tout diviseur orbifolde $\Delta$ sur $X$. Ces assertions r\'esultent du fait que pour tous entiers $m_j,t_j$, on a: $pgcd_j\{t_j\}$ divise $pgcd_j\{t_j.m_j\}$.
\end{remark}

\subsection{Groupe fondamental}\label{grofon}

On s'int\'eressera ici, lorsque $(X\vert \Delta)$ est lisse et K\" ahler, \`a son groupe fondamental $\pi_1(X\vert\Delta)$, et \`a la structure de son rev\^etement universel $\widetilde{(X\vert\Delta)}$, d\'efinis ci-dessous. 

Rappelons donc (voir [D-M 93, \S 14], [Ca 07, \S 11], ou [Cat 07, definition 4.5 p. 102]) que $\pi_1(X\vert\Delta):=\pi_1(X-Supp(\Delta))/L$, o\`u $L$ est le sous-groupe {\it normal} de $\pi_1(X-Supp(\Delta))$ engendr\'e par tous les  \'el\'ements de la forme $l_j^{m_j}$, o\`u $l_j$ est un petit lacet entourant une fois le diviseur $D_j$. (On omet les points-base). C'est donc un groupe de pr\'esentation finie puisque $\pi_1(X-Supp(\Delta))$ est de pr\'esentation finie.

\label{kertors}On a donc un morphisme naturel surjectif: $\pi_1(X\vert\Delta)\to \pi_1(X)$ dont le noyau $K$ est {\it engendr\'e par des \'el\'ements de torsion} (les images dans $\pi_1(X)$ des \'el\'ements $l_j$ ci-dessus. Voir [Ca 01], lemme 1.9.9). Cette propri\'et\'e sera cruciale dans la suite. Observer cependant que $K$ n'est en g\'en\'eral, ni de torsion, ni de type fini\footnote{Consid\'erer, par exemple, l'orbifolde hyperbolique $(E\vert \sum_{j=1}^{2g}(1-\frac{1}{2}).\{p_j\})$, o\`u $E$ est une courbe elliptique: un rev\^etement orbifolde-\'etale de degr\'e $2$ est une courbe de genre $g\geq 2$.}. Une premi\`ere cons\'equence de cette propri\'et\'e est que si $r:\pi_1(X\vert \Delta)\to G$ est une repr\'esentation dans un groupe $G$ sans torsion, elle se factorise par $\pi_1(X)$. 

Plus g\'en\'eralement, si $r:\pi_1(X\vert \Delta)\to G$ est une repr\'esentation dans un groupe r\'esiduellement fini (par exemple si $G\subset Gl(N,\bC)$ est lin\'eaire), alors sa restriction \`a un sous-groupe $G'$ d'indice fini ad\'equat a une image sans torsion, et se factorise donc par $\pi_1(X')$, o\`u $(X'\vert \Delta')$ est un rev\^etement orbifolde-\'etale fini de $(X\vert\Delta)$.

\subsection{Rev\^etement universel orbifolde}\label{runiv}

Soit $(X\vert \Delta), \Delta:=\sum_j(1-\frac{1}{m_j}).D_j$ une orbifolde g\'eom\'etrique lisse\footnote{Les \'enonc\'es qui suivent sont en fait valables, plus g\'en\'eralement, lorsque les groupes fondamentaux locaux sont finis.}. Elle admet (voir [N 87, 1.3.15, p.43]) un rev\^etement universel $r:\widetilde{(X\vert \Delta)}\to (X\vert\Delta)$, c'est-\`a-dire un rev\^etement connexe et simplement connexe, normal, Galoisien de groupe $G=\pi_1(X\vert \Delta)$, non-ramifi\'e au-dessus du compl\'ementaire de $Supp(\Delta)$, et ramifiant au-dessus de chaque point lisse $x\in D_j$ du support de $\Delta$ \`a un ordre $m^*_j$ qui divise $m_j$. Un tel rev\^etement est unique \`a automorphisme de rev\^etement pr\`es. Les entiers $m^*_j$ sont \'egalement d\'etermin\'es par ces conditions. L'orbifolde g\'eom\'etrique $(X\vert \Delta^*)$, avec $\Delta^*:=\sum_j(1-\frac{1}{m^*_j}).D_j$ est la ``r\'eduction" de $(X\vert \Delta)$ (rendue ``d\'eveloppable", par division ad\'equate de ses multiplicit\'es). On a, en particulier: $G=\pi_1(X\vert \Delta)\cong \pi_1(X\vert \Delta^*)$. Voir [N 87], [D-M 93, \S 14] ou [Ca 07, \S 11.3] pour plus de d\'etails. 

\

Par exemple, si $X=C$ est une courbe projective lisse et connexe, alors $(C\vert \Delta)$ est hyperbolique si et seulement si  son rev\^etement universel est le disque unit\'e $\Bbb D$, et $(C\vert \Delta)=\Bbb D/G$, o\`u $G$ est le $\pi_1$-orbifolde de $(C\vert\Delta)$, agissant proprement et discontinuement sur le disque unit\'e $\Bbb D$. Une telle orbifolde est ``r\'eduite".

Si $X=\Bbb P^2$ muni de deux droites projectives affect\'ees d'une m\^eme multiplicit\'e $m\geq 2$, le groupe fondamental est $\Bbb Z_m$ et le rev\^etement universel le c\^one sur la courbe rationnelle normale de degr\'e $m$, donc singulier. Cette orbifolde est ``r\' eduite".

\

On a, de plus, une bijection naturelle entre:

1. les sous-groupes $G'$ de $G:=\pi_1(X\vert \Delta)$, et:

2. Les rev\^etements {\it orbifolde-\'etales} $v:(X'\vert \Delta')\to (X\vert \Delta)$: ce sont ceux qui sont non-ramifi\'es au-dessus du compl\'ementaire de $Supp(\Delta)$, et ramifient au-dessus de chaque point lisse $x\in D_j$ du support de $\Delta$ \`a un ordre $m'_j$, avec $m_j=m'_j. m_{\Delta'}(D'_j)$, o\`u $m_{\Delta'}(D'_j)$ est la $\Delta'$-multiplicit\'e de l'unique composante $D'_j$ de $v^*(D_j)$ passant par $x'\in v^{-1}(x)$. 

\

Par exemple, si $a:(C'\vert \Delta_{C'})\to (C\vert \Delta_C)$ est un morphisme orbifolde-\'etale, alors $(C\vert \Delta_C)$ est hyperbolique si et seulement si $(C'\vert\Delta_{C'})$ l'est.

\

Un tel rev\^etement $X'$ a, au plus, des singularit\'es quotient. Si $d':Y'\to X'$ est une d\'esingularisation, elle induit donc un isomorphisme entre $\pi_1(Y')$ et $\pi_1(X')$, par [Ko 93, 7.2-7.5]. Lorsque ce rev\^etement est fini, on peut choisir $Y'$ K\" ahler, puisque $X'$ est alors dans la classe $\cal C$ de Fujiki.

Les \'eventuelles singularit\'es de l'espace normal $\widetilde{(X\vert \Delta)}$ sont donc des singularit\'es quotient (par des groupes ab\'eliens dans notre cas orbifolde lisse). Cet espace admet une d\'esingularisation $G$-\'equivariante $d:\widetilde{Y}\to \widetilde{(X\vert \Delta)}$ qui est K\" ahler si $X$ l'est. De plus, $\widetilde{(X\vert \Delta)}$ (ainsi que tous les rev\^etements orbifolde-\'etales Galoisiens interm\'ediaires) admettent une structure K\"ahl\'erienne \'equivariante au sens orbifolde usuel. Ce fait est \'etabli dans [Cl 08'].

\subsection{Morphismes sur des courbes.}

Nous aurons besoin du:

\begin{lemma}\label{orbet} Soit un diagramme commutatif:

\
\centerline{
\xymatrix{ (X'\vert\Delta' )\ar[r]^{g}\ar[d]_{h'} & (X\vert\Delta)\ar[d]^{h}\\
C'\ar[r]^a& C\\
}}

\

dans lequel les fl\`eches verticales sont des fibrations sur des courbes projectives, $X$ et $X'$ \'etant normaux et connexes compacts. On suppose que $g$ est orbifolde-\'etale. Alors, $a:(C'\vert \Delta_{h', \Delta'})\to (C\vert \Delta_{h, \Delta})$ est aussi orbifolde-\'etale. 

En particulier: $(C\vert \Delta_{h, \Delta})$ est hyperbolique si et seulement si $(C'\vert\Delta_{h', \Delta'})$ l'est.
\end{lemma}

{\bf D\'emonstration:} Soit $b\in C$, et $b'\in C'$ tel que $a(b')=b$. On note $r_a(b')$ l'ordre de ramification de $a$ en $b'$. Il s'agit de montrer que $r_a(b').m_{\Delta_{h'}}(b')=m_{\Delta_h}(b)$. 

On a: $h^*(b)=\sum_k t_k.F_k$, et $g^*(\sum_k t_k.F_k)=\sum_{k,j}t_k.r_g(F'_{k,j})+...$, en posant: $g^*(F_k)=\sum_j r_g(F'_{k,j}).F'_{k,j}+...$, en d\'esignant par $F'_{k,j}$ les composantes irr\'eductibles de $g^{-1}(F_k)$ contenues dans $(h')^{-1}(b')$, et par $r_g(F'_{k,j})$ l'ordre de ramification de $g$ le long de $F'_{k,j}$. 

Par ailleurs: $r_a(b').(h')^*(b')=r_a(b').\sum_{k,j}.t'_{k,j}.F'_{k,j}=\sum_{k,j}t_k.r_g(F'_{k,j}).F'_{k,j}$, puisque $h\circ g=a\circ h'$.

Donc: $r_a(b').t'_{k,j}=t_k.r_g(F'_{k,j}), \forall k,j$.

Puisque $g$ est orbifolde-\'etale, on a aussi: $m_{\Delta'}.r_g(F'_{k,j})=m_{\Delta}(F_k), \forall k,j$.

Par d\'efinition: $r_a(b'). m_{\Delta_h'}(b')=pgcd_{(k,j)}\{r_a(b'). t'_{k,j}.m_{\Delta'}(F'_{k,j})\}=$

=$pgcd_{(k,j)}\{t_k.r_g(F'_{k,j}).m_{\Delta'}(F'_{k,j})\}=pgcd_{(k,j)}\{t_k.m_{\Delta}(F_{k})\}=m_{\Delta_h}(b)$ $\square$

\

Signalons enfin que si $f:X\to C$ est un morphisme surjectif de $X$, compacte K\"ahler, sur une courbe projective $C$, et si $(X\vert\Delta)$ est lisse, alors $f$ induit une suite exacte de groupes: $$1\to \pi_1(F\vert \Delta_F)\to G\to \pi_1(C\vert \Delta_{f,\Delta})\to 1,$$ 
dans laquelle $\Delta_F$ est la restriction de $\Delta$ \`a la fibre g\'en\'erique $F$ de $f$. Voir [Ca 07, Prop.11.7] pour la d\'emonstration (en toute dimension pour $C$, pourvu que le morphisme $f$ soit `net',  ce qui peut \^etre r\'ealis\'e par modifications ad\'equates de $X$ et $C)$.

\

\section{Ensembles de Green-Lazarsfeld des orbifoldes g\'eom\'etriques.}

\subsection{Ensembles de Green-Lazarsfeld: g\'en\'eralit\'es.}\label{gl}

Soit $G$ un groupe de type fini, $DG$ son groupe d\'eriv\'e (engendr\'e par ses commutateurs), $G_{ab}=G/DG$ son ab\'elianis\'e. 

Le groupe des caract\`eres complexes (multiplicatifs) de $G$ est $\hat G:=Hom(G,\bC^*)=H^1(X,\Bbb C^*)$, par le th\'eor\`eme d'Hurwicz. On notera $\chi$ un tel caract\`ere, et $\bC_{\chi}$ le $G$-module dont le groupe additif est $\bC$, muni de l'action: $g.u:=\chi(g).u$, pour $(u,g)\in \bC\times G$.

On note $\Sigma^1(G)\subset \hat G$ {\it l'ensemble de Green-Lazarsfeld de $G$}, constitu\'e des caract\`eres $\chi:G_{ab}\to \bC^*$ tels que $H^1(G,\bC_{\chi})\neq 0$. La suite spectrale de Hochschild-Serre fournit un isomorphisme naturel $H^1(G,\bC_{\chi})\cong Hom_G((DG)_{ab},\bC_{\chi})$, $G$ agissant sur $(DG)_{ab}$ par conjugaison, et sur $\bC$ par multiplication par $\chi$ ([Be 92, 4.2]).

Un morphisme de groupes surjectif $q:G\to Q$ induit des injections $q^*:\hat Q\to \hat G$ et $\Sigma^1(Q)\to \Sigma^1(G)$.

\subsection{Ensembles de Green-Lazarsfeld orbifoldes.}\label{glorb'}

Si $(X\vert \Delta)$ est une orbifolde g\'eom\'etrique finie, lisse et (compacte) K\" ahler, on notera $G:=\pi_1(X\vert \Delta)$, $H:=\pi_1(X)$, et $f:G\to H$ le morphisme de groupes surjectif naturel, de noyau $K$. On posera: $\Sigma^1(X\vert\Delta):=\Sigma^1(G)\subset \hat G=H^1(X,\bC^*)$. On a une inclusion \'evidente $f^*:\hat H\to \hat G$ dont le conoyau est fini, puisque $K$ est engendr\'e par des \'el\'ements de torsion, et que $G$ est de type fini.

L'ensemble $\Sigma^1(G)$ ne d\'epend que du m\'etab\'elianis\'e de $G$, c'est-\`a-dire du quotient de $G$ par son second groupe d\'eriv\'e $D^2(G)=D(D(G))$. Son int\'er\^et dans le cas o\`u $G=\pi_1(X\vert \Delta)$ est que sa structure permet de montrer que ce m\'etab\'elianis\'e est virtuellement nilpotent si $DG$ est de type fini, et que $(X\vert \Delta)$ fibre sur une courbe orbifolde hyperbolique sinon.

La th\'eorie de Hodge fournit de cruciales restrictions sur sa structure (d\'ecouvertes dans [GL 87]) d\'ecrites ci-dessous, qui impliquent de tr\`es fortes restrictions sur la structure du groupe $G$ lui-m\^eme (observ\'ees pour la premi\`ere fois dans [Be 92]). Ces implications proviennent de la suite exacte:
$$1\to DG/D^2G\to G/D^2G\to G/DG:=G_{ab}\to 1,$$
et de l'action par conjugaison de $G_{ab}$ sur $DG/D^2G$. Pour cette action, les composants isotypiques sont pr\'ecis\'ement les \'el\'ements de $\Sigma^1(G)$. Ces observations seront d\'evelopp\'ees ci-dessous.

\begin{theorem}\label{glorb} Soit $(X\vert \Delta)$ une orbifolde g\'eom\'etrique finie, lisse et K\" ahler, et $\Sigma^1(X\vert \Delta)\subset H^1(X,\bC^*)$ son ensemble de Green-Lazarsfeld. Alors:

1. $\Sigma^1(X\vert \Delta)$ est une r\'eunion finie de translat\'es de sous-tores par des \'el\'ements de torsion.

2. il existe un nombre fini de morphismes orbifoldes $f_j$ de $(X\vert \Delta)$ sur des courbes orbifoldes hyperboliques $(C_j\vert \Delta_{C_j})$ de genre $g(C_j)>0$ tels que la r\'eunion des composantes irr\'eductibles de dimension strictement positive de $\Sigma^1(X\vert \Delta)$ soit la r\'eunion des $(f_j)^*(\Sigma^1(C_j\vert\Delta_{C_j}))$.

\end{theorem}

Si $DG$ est de type fini, il n'existe pas (voir 4.5.(2)) de morphisme $f_j$ comme ci-dessus, donc:

\begin{corollary}\label{gl'} Dans la situation du th\'eor\`eme \ref{glorb}, supposons que $DG$ soit de type fini, avec $G:=\pi_1(X\vert \Delta)$. Alors: $\Sigma^1(X\vert \Delta)$ est un ensemble fini compos\'e d' \'el\'ements de torsion.
\end{corollary}

\begin{remark} La situation est donc l'exact analogue du cas o\`u $\Delta=0$ trait\'e dans [GL 87], [Be 92], [Si 93] et [Ca 01]. La d\'emonstration donn\'ee ci-dessous consiste \`a se ramener \`a ce cas particulier par des arguments dont certains sont analogues \`a ceux utilis\'es dans [Ca 01, (1.9.3-1.11.8, pp. 603-607) ] pour d\'eduire le cas K\" ahler du cas projectif. Dans le \S 3.3, nous donnons aussi, en suivant une suggestion du rapporteur, l'esquisse d'une preuve directe, enti\`erement ind\'ependante, qui consiste \`a adapter au cadre orbifolde la d\'emonstration de [De 07]. 

Remarquons aussi que l'on pourrait proc\'eder ci-dessous de mani\`ere plus naturelle et directe, en r\'eecrivant dans le cadre des orbifoldes g\'eom\'etriques (ou de la th\'eorie des champs, ou des V-vari\'et\'es de Satake) la d\'emonstration des r\'esultats de [Be 92] et de [G-L 87]. Nous y avons renonc\'e: faute de r\'ef\'erences appropri\'ees, il aurait \'et\'e n\'ecessaire de red\'efinir les notions utilis\'ees dans ce cadre, ce qui \'etait en soi d\'ej\`a un article ind\'ependant.

\end{remark}

{\bf D\'emonstration du th\'eor\`eme \ref{glorb}:} Elle consiste \`a construire un rev\^etement orbifolde-\'etale fini $g:(X'\vert \Delta')\to (X\vert \Delta)$ tel que la restriction des \'el\'ements de $\Sigma^1(G)$ proviennent d'\'el\'ements de $\Sigma^1(X')$ (d\'esingularis\'ee), \`a exploiter la structure (connue) de ce dernier, puis de redescendre ces constructions sur $(X\vert\Delta)$.

Plus pr\'ecis\'ement, nous allons construire un diagramme commutatif:

\
\centerline{
\xymatrix{ Y'\ar[rd]^{d'}&(X'\vert \Delta')\ar[r]^{g}\ar[d]_{h}&(X\vert \Delta)\ar[dd]_{f}\\&X'\ar[d]_{f'} &\\
&C'\ar[r]^a&C\\
}}

tel que:

1. $g$ soit orbifolde-\'etale, et qu'existent $\chi'\in \Sigma^1(X')$ et $u'\in H^1(X',\bC_{\chi'})$ tels que $g^*(\chi)=h^*(\chi')$, et $g^*(u)=h^*(u')$ pour tout $\chi\in \Sigma^1(G)$ et $u\in H^1(G,\chi)$. En particulier, $g^*(\Sigma^1(G)\subset h^*(\Sigma^1(X'))$. C'est la premi\`ere \'etape (proposition \ref{p1}). 

2. Nous distinguons alors deux cas: $\chi$ est isol\'e dans $\Sigma^1(G)$ et donc de torsion, ou bien $\chi$ n'est pas isol\'e dans $\Sigma^1(G)$. Dans ce cas, que l'on suppose r\'ealis\'e dans la suite, on construit les morphismes $f',a,f$ du diagramme ci-dessus, et on montre que les bases-orbifoldes de $f'\circ h$ et de $f$ sont hyperboliques. C'est la seconde \'etape (lemmes \ref{C} et \ref{chyp}). 

3. Dans la troisi\`eme \'etape, on montre (lemmes \ref{desc} et \ref{glcourb}) l'existence de $\gamma\in \Sigma^1(\pi_1(C\vert \Delta_{f,\Delta}))$ tel que $\chi=f^*(\gamma)$, ce qui ach\`eve la d\'emonstration.

\

 La premi\`ere \'etape est la suivante:

\begin{proposition}\label{p1}  Il existe un rev\^etement orbifolde-\'etale fini $g:(X'\vert \Delta')\to (X\vert \Delta)$ tel que, pour tout couple $(\chi,u)$ avec $1\neq\chi\in \Sigma^1(G)$ et $0\neq u\in Hom(DG_{ab},\Bbb C_{\chi})$, si $H':=\pi_1(X')$, alors il existe $\chi'\in \Sigma^1(H')$ et $0\neq u'\in Hom(DH'_{ab},\Bbb C_{\chi'})$ tels que $g^*(\chi)=h^*(\chi')$ et $g^*(u)=h^*(u')$, si $h:G':=\pi_1(X'\vert\Delta')\to H'$ est le morphisme de groupes naturel.
\end{proposition}

{\bf D\'emonstration:} Posons $G:=\pi_1(X\vert \Delta)$, et $H:=\pi_1(X)$. Le morphisme de groupes naturel surjectif $f:G\to H$ a un noyau $K$ engendr\'e par des \'el\'ements de torsion (voir \ref{kertors}). Il induit donc un morphisme de groupes (ab\'eliens) surjectif $(Df)_{ab}:(DG)_{ab}\to (DH)_{ab}$.

\begin{lemma}\label{revet} (Voir [Ca 01, 1.9.8]) Avec les notations pr\'ec\'edentes, on a:

1.  $DG$ est d'indice fini dans $K.DG$.

2. Il existe un sous-groupe $G'\subset G$ normal et d'indice fini, avec $G/G'$ ab\'elien, tel que $K\cap G'\subset DG$.
\end{lemma}

{\bf D\'emonstration de \ref{revet}:} L'application compos\'ee $K\to G\to G_{ab}$ a pour noyau $K\cap DG$, et une image $T$, isomorphe \`a $K.DG/DG$, qui est finie, puisque $K$ est engendr\'e par des \'el\'ements de torsion, et $G_{ab}$ ab\'elien de type fini. D'o\`u l'assertion 1. 

Soit $G"\subset G_{ab}$ un suppl\'ementaire dans $G_{ab}$ (donc sans torsion) du sous-groupe de torsion $Tors(G_{ab})\supset T$. C'est-\`a-dire que l'on a: $G_{ab}=G"\oplus Tors(G_{ab})$, et en particulier: $G"\cap T=\{1\}$. On d\'efinit alors $G'$ comme \'etant l'image r\'eciproque de $G"$ dans $G$. Alors: $K\cap G\subset DG$. En effet: l'image de $G'$ dans $G_{ab}$ est sans torsion, tandis que celle de $K$ est finie $\square$

\

 Nous choisissons et fixons maintenant un sous-groupe normal et d'indice fini $G'\subset G$ comme dans le lemme pr\'ec\'edent. On note $g^*:\widehat{G}\to\widehat{G'}$ la restriction naturelle (de noyau fini).

\begin{lemma}\label{indfin} $g^*(\Sigma^1(G))\subset \Sigma^1(G')$, et $g^*:H^1(G,\chi)\to H^1(G',g^*(\chi))$ est injective, pour tout $\chi\in \hat G$.
\end{lemma}

{\bf D\'emonstration} La corestriction (voir [Se 68, VII,\S 7, proposition 6, p. 127]) montre cette injectivit\'e, puisque $G'$ est d'indice fini dans $G$. La premi\`ere assertion en r\'esulte. $\square$ 





\

 Soit $g:(X'\vert \Delta')\to (X\vert \Delta)$ le rev\^etement {\it orbifolde-\'etale} fini d\'efini par l'inclusion $G'\subset G$ de \ref{revet} (voir \ref{runiv} ci-dessus). Soit $h':G'\to H':=\pi_1(X')$ le morphisme de groupes naturel. Son noyau $K'$ est donc engendr\'e par des \'el\'ements de torsion. De plus, $K'\subset (K\cap G')\subset DG$. 
 
 Par le lemme \ref{indfin} ci-dessus, si $\chi\in \widehat{G}$, et si $0\neq u\in Hom_G(DG_{ab},\bC_{\chi})$, alors $0\neq g^*(u)\in Hom_{G'}((DG')_{ab},\bC_{g^*(\chi)})$. 
 
 Par le lemme \ref{facto} ci-dessous (appliqu\'e \`a $G'$, $N=K'$, et $g^*(\chi)$, avec $G/N=H')$, on a:$g^*(\chi)=(h')^*(\chi')$, pour $\chi'\in \widehat{H'}$, et $g^*(u)=(h')^*(u')$, pour $0\neq u'\in H^1(H',\bC_{\chi'})$. Donc $\chi'\in \Sigma^1(X')$. Ce qui ach\`eve la d\'emonstration de la proposition \ref{p1} $\square$

 \

 \begin{lemma}\label{facto} Soit $\chi\in \Sigma^1(G)$, et $u\in Hom_G((DG)_{ab},\bC_{\chi})$. Soit $N\subset DG$ un sous-groupe normal de $DG$ tel que $u$ s'annule sur $N$. Alors:

1. Il existe $\chi'\in H^1(G/N,\bC)$ tel que $\chi=q^*(\chi')$, si $q: G\to Q:=G/N$ est le quotient.

2. Il existe $v\in Hom_Q((DQ)_{ab}, \bC_{\chi'})$ tel que $u=q^*(v)$ (dans un sens fonctoriel \'evident)
\end{lemma}

{\bf D\'emonstration de \ref{facto}:} Il suffit de v\'erifier l'existence de $\chi'$, c'est-\`a-dire que $\chi(k)=1$ pour $k\in N$, ce qui est \'evident, puisque  $\chi(k)=1$ pour $k\in DG$, et que $N\subset DG$. $\square$

\

 La proposition \ref{p1}, et la premi\`ere \'etape, sont ainsi \'etablies.

 \begin{remark}\label{d} Si $d':Y'\to X'$ est une d\'esingularisation K\" ahl\'erienne, elle induit, d'apr\`es [Ko 93,7.2-7.5], un isomorphisme de groupes au niveau des groupes fondamentaux. Les ensembles de Green-Lazarsfeld de $X'$ et $Y'$ seront donc identifi\'es ici. Autrement dit, on peut supposer $X'$ lisse, et $\Sigma^1(X')$ a donc la structure d\'ecrite dans l'\'enonc\'e \ref{glorb}, par les travaux cit\'es ci-dessus. 

Remarquons aussi que, puisque les singularit\'es de $X'$ sont quotient, si $f':Y'\to C'$ est un morphisme sur une courbe $C'$ de genre $g'>0$, alors ce morphisme se factorise par $d'$.
\end{remark}

 \

 Nous abordons maintenant la seconde \'etape de la d\'emonstration du th\'eor\`eme \ref{glorb}. 
   Distinguons deux cas:
 
 \
 
 {\bf 1.} $\chi'$ est isol\'e dans $\Sigma^1(X')$, et donc de torsion, par [Ca 01] (qui r\'eduit le cas K\" ahler au cas projectif, \'etabli par C. Simpson dans [Si 93]). Alors $\chi$ est de torsion dans $\widehat G$.  
 
 {\bf 2.} $\chi'$ n'est pas isol\'e dans $\Sigma^1(X')$. 
 C'est ce second cas qu'il nous reste donc \`a traiter.
 
 \

 
 Par [Be 92, 3.5] reformul\'e dans la terminologie orbifolde, et les remarques \ref{d} et \ref{d'} ci-dessus, pour toute composante irr\'eductible de $\Sigma^1(X')$ contenant $h^*(\chi)$, il existe un morphisme orbifolde $f':X'\to (C'\vert \Delta_{f'})$ tel que la composante irr\'eductible consid\'er\'ee de $\Sigma^1(X')$ contenant $h^*(\chi')$ soit \'egale \`a une composante irr\'eductible de $(f')^*(\Sigma^1(\pi_1(C'\vert \Delta_{f'}))$. 
 Fixons une telle composante.

 \

 Nous allons d'abord montrer que $f'$ provient d'un morphisme $F:X\to C$.
 
\begin{lemma}\label{C} Il existe un diagramme commutatif, avec $C$ une courbe, $F$ et $a$ \'etant des morphismes surjectifs:
 
 \
\centerline{
\xymatrix{ X'\ar[r]^{g}\ar[d]_{f'} & X\ar[d]^{F}\\
C'\ar[r]^a& C\\
}}

\end{lemma}

{\bf D\'emonstration:} Consid\'erons le diagramme commutatif:

\
\centerline{
\xymatrix{ B'\ar[r]^{g}\ar[d]_{}&B\ar[d]_{}\\Alb(X')\ar[r]^{g}\ar[d]_{f'} & Alb(X)\ar[d]^{f}\\
Alb(C')\ar[r]^a&Alb(X)/B\\
}}

\

Dans lequel on utilise les m\^emes lettres pour d\'esigner une application, l'une de ses restrictions,  et l'application induite au niveau des vari\'et\'es d'Albanese. De plus, $B'$ d\'esigne le noyau (qui est connexe) du morphisme $f'$, et $B:=g(B')$.

Dans ce diagramme, $B\neq Alb(X)$. En effet: $\chi'=(f')^*(\gamma')$ est trivial sur le noyau $\pi_1(B')$, et $\chi$ est donc trivial sur $\pi_1(B)$, puisque $g^*(\chi)=\chi'$. Si $B=Alb(X)$, $\chi$ est donc de torsion, contrairement \`a notre hypoth\`ese. (Remarquer que $\widehat{\pi_1(X)}$ et $\widehat{\pi_1(Alb(X))}$ ont la m\^eme composante neutre).

Donc: $C:=a(C')$ est une courbe irr\'eductible de $Alb(X)/B$, puisque $C'$ engendre $Alb(C')$, de sorte que $C$ engendre $Alb(X)/B\neq \{0\}$. La commutativit\'e du diagramme pr\'ec\'edent montre que $f:X\to C$ est surjective, avec $f\circ g=a\circ f':X'\to C$. (On a normalis\'e $C$, et encore not\'e $a$ la restriction de $a$ \`a $C'$, enfin $\alpha_X$ est le morphisme d'Albanese de $X$)
$\square$

\begin{lemma}\label{chyp} La base orbifolde de $F:(X\vert \Delta)\to C$ est hyperbolique.
\end{lemma}

{\bf D\'emonstration:} On d\'eduit des morphismes $f':X'\to C'$ et $f:X\to C$ le diagramme commutatif suivant de morphismes orbifoldes, dans lequel $g$ est, par construction, orbifolde-\'etale:

\
\centerline{
\xymatrix{ (X'\vert \Delta')\ar[r]^{g}\ar[d]_{F'}&(X\vert \Delta)\ar[d]_{F}\\(C'\vert\Delta_{f',\Delta'})\ar[r]^{a}& (C\vert \Delta_{F,\Delta})\\
}}

\

On d\'eduit alors du lemme \ref{orbet} que $a:(C'\vert\Delta_{f',\Delta'})\to (C\vert \Delta_{F,\Delta})$ est orbifolde-\'etale, et donc que $(C\vert \Delta_{F,\Delta})$ est une courbe orbifolde hyperbolique si $(C'\vert\Delta_{f',\Delta'})$ l'est. V\'erifions cette derni\`ere propri\'et\'e.

Soit $F":=f'\circ d':Y'\to C'$; nous savons que $(C'\vert \Delta_{F"})$ est hyperbolique. Donc $(C'\vert \Delta_{f'})$ est {\it a fortiori} hyperbolique, puisque $\Delta_{F"}$ divise \'evidemment $\Delta_{f'}$(voir la derni\`ere remarque de la section \ref{orbgeo}). Enfin, $\Delta_{f'}$ divise $\Delta_{f',\Delta'}$ (\'egalement par \ref{orbgeo}), ce qui \'etablit l'hyperbolicit\'e de $(C'\vert \Delta_{f',\Delta'})$. $\square$

\

 Troisi\`eme \'etape.

 Rappelons que $g^*(\chi)=h^*(\chi')$, pour $\chi'\in \Sigma^1(H')$, $H':=\pi_1(X')$, et que $\chi'$ n'est pas isol\'e dans $\Sigma^1(H')$. On a, de plus, identifi\'e $H'$ et $\pi_1(Y')$ si $d':Y'\to X'$ est une d\'esingularisation K\" ahl\'erienne. On note aussi: $F":=f'\circ d':Y'\to C'$.

\

Le th\'eor\`eme \ref{glorb} r\'esulte alors des lemmes \ref{desc} et \ref{glcourb} ci-dessous:

\begin{lemma}\label{desc} Posons $M:=\pi_1(C\vert\Delta_{F,\Delta})$, et $M':=\pi_1(C'\vert\Delta_{F'})$.

Si $\chi'\in \Sigma^1(H')$ n'est pas isol\'e dans $\Sigma^1(H')$, alors:

1.  $(d')^*\circ g^*(\chi)=(f')^*(\gamma'),$ pour un $ \gamma'\in\Sigma^1(M').$

2. $H^1(H',\bC_{\chi'})=(f')^*(H^1(M',\bC_{\gamma'})).$

3. il existe $\gamma\in \Sigma^1(M)$ tel que $\chi=F^*(\gamma)$, et $H^1(G,\bC_{\chi})\subset F^*(H^1(M,\bC_{\gamma}))$.
\end{lemma}

{\bf D\'emonstration:} Assertion 1. Nous avons vu que $(d')^*\circ g^*(\chi)=h^*(\chi')$ pour un $\chi'\in \Sigma^1(H')$. Par [Be 92], la composante connexe de $\Sigma^1(H')$ contenant $\chi'$, suppos\'e non isol\'e dans $\Sigma^1(H')$ est une composante connexe de $(f')^*(\Sigma^1(M'))$. D'o\`u l'assertion, puisque $(F')^*:=(d')^*\circ (f')^*$, et que $d'$ identifie $H'$ et $\pi_1(Y')$, et donc $\chi'$ et $(d')^*(\chi')$.

\

Assertion 2. L'assertion r\'esulte de [A 92, theorem 3, p. 311].

Nous allons aussi en donner une d\'emonstration utilisant plusieurs des arguments de [Be 92]. Nous avons tout d'abord l'\'egalit\'e: $ H^1(H',\bC_{\chi'})= H^1(Y',\cal L')$, $\cal L'$ d\'esignant le faisceau localement constant sur $Y'$ d\'eduit de $\bC_{\chi'}$, ainsi qu'une suite exacte courte (voir [Be 92, d\'emonstration de la prop.3.5], et [A 92, Theorem 3]) :
$$0\to E_2^{1,0}(Y',\chi')\to H^1(H',\bC_{\chi'})\to E_2^{0,1}(Y',\chi')\to 0,$$
 dans laquelle $E_2^{1,0}(Y',\chi')$ est l'homologie du complexe: $$H^0(Y',L')\to H^0(Y',\Omega_{Y'}^{1}\otimes L')\to H^0(Y',\Omega_{Y'}^{2}\otimes L'),$$
 tandis que $E_2^{0,1}(Y',\chi')$ est le noyau de la fl\`eche: $$H^1(Y',L')\to H^1(Y',\Omega_{Y'}^{1}\otimes L'),$$ dans lesquels les fl\`eches sont donn\'ees par $\wedge \omega$, tandis que, si $\cal L'$ est le fibr\'e en droites plat sur $Y'$ donn\'e par le caract\`ere $\chi'$, dont les sections locales sont constantes, alors $L'$ d\'esigne le fibr\'e en droites holomorphe d\'efini par: $L':=\cal L'$$\otimes_{\bC}\cal O$$_{Y'}$. Enfin, $\omega'$ est une $1$-forme holomorphe non-nulle sur $C'$.
 
 \

 Par fonctorialit\'e, nous avons un diagramme commutatif:
 
 \
 
 \
\centerline{
\xymatrix{0\ar[r]& E_2^{1,0}(C',\gamma')\ar[r]^{\wedge \omega}\ar[d]_{F'^*}& H^1(M',\bC_{\gamma'})\ar[r]^{\wedge \omega}\ar[d]_{F'^*}& E_2^{0,1}(C',\gamma')\ar[r]\ar[d]_{F'^*}&0\\
0\ar[r]^{}& E_2^{1,0}(Y',\chi')\ar[r]^{\wedge F'^*(\omega)}& H^1(H',\bC_{\chi'})\ar[r]^{\wedge F'^*(\omega)}& E_2^{0,1}(Y',\chi')\ar[r]& 0\\
}}

\

\

Nous voulons montrer que la fl\`eche verticale du milieu est bijective. Il suffit de montrer que les deux autres fl\`eches verticales extr\^emes le sont. Pour celle de gauche, ceci r\'esulte des arguments donn\'es dans la d\'emonstration de [Be 92, Prop. 2.1]. Pour la fl\`eche de droite, ceci r\'esulte du fait que $E_2^{0,1}(Y',\chi')$ est le conjugu\'e de $E_2^{1,0}(Y',(\chi')^{-1})$, ce qui est d\'emontr\'e (mais non formul\'e) dans la preuve de [Be 92, Prop. 3.5].

 
\

Assertion 3. Soit en effet $0\neq u\in Hom_G(DG,\bC_{\chi})$, alors $g^*(u)$ ne s'annule pas sur $DG'$par le lemme \ref{indfin}. 
Puisque, par l'assertion 2 pr\'ec\'edente, $g^*(u)=(F')^*(v)$, $v\in H^1(C',\bC_{\gamma'})$, $u$ s'annule sur $\Phi':=\pi_1(\Psi'), \Psi'$ \'etant une fibre lisse de $f':X'\to C'$. Donc $u$ s'annule sur $g_*(\Phi')$, qui est un sous-groupe d'indice fini de $\Phi:=\pi_1(\Psi)$, $\Psi:=g(\Psi')$ \'etant une fibre lisse de $f:X\to C$. 

Puisque le groupe (additif) $\bC_{\chi}$ est sans torsion, $u$ s'annule aussi sur $\Phi$, et m\^eme sur $\Phi_{\Delta}:=\pi_1(\Psi\vert \Delta_{\Psi})$. On a not\' e $\Delta_{\Psi}:= \Delta_{\vert \Psi}$ la restriction du diviseur $\Delta$ \`a une fibre g\'en\'erique $\Psi$ de $f$.  La derni\`ere assertion provient de ce que le noyau du morphisme de groupes naturel $\Phi_{\Delta}\to \Phi$ est engendr\'e par des \'el\'ements de torsion.

On dispose enfin ([Ca 07, prop.11.7]) d'une suite exacte de groupes (car tout morphisme sur une courbe est {\it net}) induite par $f_*$:

$$\Phi_{\Delta}\to G\to M\to  1$$

Le lemme \ref{facto} montre alors que $\chi=F^*(\gamma)$, pour $\gamma\in \widehat{M}$. La seconde assertion de 3 est \'evidente $\square$

\begin{lemma}\label{glcourb} Soit $C$ une courbe projective lisse et connexe de genre $g(C)\geq 1$. Soit $\Delta_C$ une structure orbifolde sur $C$ telle que $(C\vert \Delta_C)$ soit hyperbolique.  Soit $M:=\pi_1(C\vert \Delta_C)$. Alors: $\Sigma^1(M)=\widehat{M}$ si $g(C)\geq 2$, et $\Sigma^1(M)=\widehat{M}-\widehat{M}^0$ si $g(C)=1$. On d\'esigne ici par $\widehat{M}^0=H^1(C,\bC^*)$ la composante neutre de $\widehat{M}$, priv\'ee de son \'el\'ement unit\'e.
\end{lemma}

{\bf D\'emonstration:} Cet \'enonc\'e r\'esulte des arguments des d\'emonstrations de [Be 92, \S 1 et Prop. 3.5] $\square$

\begin{remark}\label{rbe} Dans la situation de \ref{glcourb}, si $\Delta_C=\sum_j(1-\frac{1}{m_j}).\{p_j\}$, le groupe de torsion $T$ de $\widehat{M}$ est isomorphe au sous-groupe $T(\Delta)$ de $\oplus_j\bZ/m_j$ constitu\'e des familles de classes $(\bar {k_j })_j$ telles que $\sum_j \frac{k_j}{m_j}$ soit entier. (Voir [Be 92, 1.5]).
\end{remark}

\subsection{Seconde d\'emonstration, ind\'ependante.}

\

Nous donnons ici, suivant une suggestion due \`a un rapporteur, une seconde d\'emonstration potentielle du th\'eor\`eme 3.1 pr\'ec\'edent. Cette d\'emonstration est une simple adaptation au cadre orbifolde de celle de [De 07]. Observons qu'elle est ind\'ependante de [Si 93] en particulier. 

Soit $\chi \in \Sigma^1(X\vert \Delta)$. Soit $K$ le sous-corps de $\Bbb C$ engendr\' e par $\chi(G), G:=\pi_1(X\vert \Delta)$, et $\sO_K$ l'anneau des entiers alg\'ebriques de $K$. Si $\chi$ n'est pas de torsion, alors le lemme de Kronecker montre que, ou bien:

1. $\chi(G)$ n'est pas contenu dans $\sO_K$, ou bien:

2. $\vert\chi(G)\vert\neq 1$.

Il suffit donc de montrer que, dans chacun de ces deux cas, $\chi$ provient d'un morphisme orbifolde sur une courbe orbifolde hyperbolique.

Dans le cas 2, ceci r\'esulte en particulier de la d\'emonstration donn\'ee ci-dessus du th\'eor\`eme 3.1 (plus pr\'ecis\'ement, de sa r\'eduction au cas o\`u $\Delta=0)$. Une d\'emonstration plus naturelle pourrait \^etre obtenue en formulant les arguments de [Be 92] dans le cadre orbifolde. (Ceci n\'ecessiterait de d\'evelopper et introduire les notions et propri\'et\'es de base utilis\'ees). 

Dans le cas 1, il est montr\'e dans [De 07] qu'il existe une valuation discr\`ete $\nu$ sur $K$ telle que le morphisme compos\'e $\nu\circ \chi: G\to \Bbb Z$ soit `exceptionnel' dans le sens de Bieri-Neumann-Strebel (c'est-\`a-dire tel que que l'image r\'eciproque dans le graphe de Cayley de $G$ (pour un syst\`eme fini arbitraire de g\'en\'erateurs) d'un intervalle $]t,+\infty]$ ad\'equat soit non-connexe). L'existence de $\nu$ repose sur un lemme alg\'ebrique utilisant le fait que l'ab\'elianis\'e de $G$ est de type fini, et donc que son alg\`ebre de groupe \`a coefficients rationnels est noeth\'erienne. Le fait que la classe de la compos\'ee reste dans $\Sigma^1(G)$ se v\'erifie en faisant agir $G$ sur l'arbre de Bruhat-Tits de $Sl(2,K)$ par l'interm\'ediaire du caract\`ere $\chi$ et de la classe de $0\neq u\in H^1(G,\chi)$. 

La d\'emonstration du th\'eor\`eme 4.5 ci-dessous montre alors que $\nu\circ \chi$ provient encore d'une courbe orbifolde hyperbolique. La d\'emonstration de la proposition 2 de [De 07] montre que $\chi$ lui-m\^eme provient de cette courbe orbifolde. 

Nous renvoyons \`a [De 07] pour les d\'etails omis.

\

\section{Quotients r\'esolubles des groupes de K\" ahler orbifoldes.}\label{quores}

Nous suivrons, en les adaptant au contexte orbifolde, les arguments de [De 06]. Le r\'esultat principal est le suivant:

\begin{theorem}\label{resnilp} Soit $(X\vert\Delta)$ une orbifolde g\'eom\'etrique lisse, avec $X$ K\"ahler compacte et connexe. On a \'equivalence entre les propri\'et\'es suivantes:

1. Tout quotient r\'esoluble $R$ de tout sous-groupe d'indice fini $G$ de $\pi_1(X\vert \Delta)$ est virtuellement nilpotent (ie: admet un sous-groupe d'indice fini nilpotent).

2. Il n'existe pas de morphisme orbifolde {\bf divisible} $f':(X'\vert\Delta')\to (C\vert\Delta_C)$ sur une courbe orbifolde hyperbolique, si $u:(X'\vert \Delta')\to (X\vert \Delta)$ est un morphisme orbifolde \'etale fini arbitraire. 

3. Il n'existe pas de morphisme de groupes surjectif: $\pi_1(X'\vert \Delta')\to \pi_1(C\vert \Delta)$, avec $(C\vert \Delta)$ une courbe orbifolde hyperbolique, si $u:(X'\vert \Delta')\to (X\vert \Delta)$ est un morphisme orbifolde \'etale fini arbitraire. 
\end{theorem}

\begin{remark} Il n'est pas possible de simplifier les \'enonc\'es ci-dessus en omettant de consid\'erer les rev\^etements orbifolde-\'etales finis de $(X\vert \Delta)$. Voir l'exemple \ref{descente}.
\end{remark}

Les implications $1\Longrightarrow3\Longrightarrow 2$ sont \'evidentes. L'implication $2\Longrightarrow 1$ r\'esulte des trois \'enonc\'es \ref{metab}, \ref{dgtf} et \ref{dgti} suivants:

\begin{theorem}\label{metab}([De 06, th\'eor\`eme 3.2]) Soit $R$ un groupe r\'esoluble de type fini non virtuellement nilpotent. Il existe alors un sous-groupe $R'$ d'indice fini de $R$ ayant un quotient m\'etab\'elien $Q'$ non virtuellement nilpotent. (Un groupe $Q'$ est m\'etab\'elien si son groupe d\'eriv\'e est ab\'elien, ie: si $D(DG):=D^2G=\{1\})$.\end{theorem}

\begin{theorem}\label{dgtf} Soit $(X\vert\Delta)$ une orbifolde g\'eom\'etrique lisse, avec $X$ K\"ahler compacte et connexe. Si le groupe d\'eriv\'e $DG$ de $G:=\pi_1(X\vert\Delta)$ est de type fini, tout quotient r\'esoluble de $\pi_1(X\vert\Delta)$ est virtuellement nilpotent.
\end{theorem}

{\bf D\'emonstration (de \ref{dgtf}):} Il nous suffit, d'apr\`es \ref{metab}, de montrer que $G/D^2G$ est virtuellement nilpotent (avec $G:=\pi_1(X\vert\Delta)$). Or $G$ agit par conjugaison de mani\`ere quasi-unipotente sur $DG/D^2G$. En effet, les valeurs propres de l'action de $G$ sur $(DG_{ab})\otimes \bC$, qui est de dimension complexe finie, ne sont autres que les \'el\'ements de $\Sigma^1(X)$.  Or ces \'el\'ements sont en nombre fini, et de torsion, d'apr\`es \ref{gl'}. Si $G'\subset G$ est le noyau de l'action de $G$ sur $(DG_{ab})$ (qui est ab\'elien de type fini), $G'$ est donc d'indice fini dans dans $G$, et $G'$ agit trivialement sur $(G'/D^2G)$, qui est donc nilpotent et d'indice fini dans $G/D^2G$ $\square$

\begin{theorem}\label{dgti} Soit $(X\vert\Delta)$ une orbifolde g\'eom\'etrique lisse, avec $X$ K\"ahler compacte et connexe. On a \'equivalence entre les propri\'et\'es suivantes:

1. Il existe $G'<\pi_1(X\vert \Delta):=G$ d'indice fini tel que le groupe d\'eriv\'e $DG'$ {\bf n'est pas} de type fini.

2. Il existe un morphisme orbifolde \'etale fini $u:(X'\vert \Delta')\to (X\vert \Delta)$, et un morphisme orbifolde {\bf divisible} $f':(X'\vert\Delta')\to (C\vert\Delta_C)$ sur une courbe orbifolde hyperbolique.

\end{theorem}


{\bf D\'emonstration:} L'implication $2\Longrightarrow 1$ est claire, puisque $D\pi_1(C\vert \Delta_C)$ n'est pas de type fini si $(C\vert \Delta_C)$ est hyperbolique et de genre $g>0$. Lorsque $g=0$, un rev\^etement orbifolde-\'etale fini de $(C\vert \Delta_C)$, suppos\'ee hyperbolique, est de genre strictement positif.

L'implication oppos\'ee est, pour l'essentiel, une adaptation au cas orbifolde de celle de [De 06, th\'eor\`eme 1.1], dont nous adopterons les notations. On peut \'evidemment supposer (quitte \`a remplacer $(X\vert \Delta)$ par un rev\^etement orbifolde-\'etale fini ad\'equat) que $D\pi_1(X\vert \Delta)$ n'est pas de type fini. Il existe donc alors (dans la terminologie de [De 06]) une classe ``exceptionnelle" r\'eelle (ie: \`a valeurs r\'eelles) $\chi\in \widehat G$, avec $G:=\pi_1(X\vert \Delta)$. Dire que cette classe est ``exceptionnelle" signifie, par d\'efinition (voir[B-S]), que $\chi^{-1}([0,+\infty[)$ n'est pas connexe (ie: dans le graphe de Cayley $C(G)$ de $G$, le sous-graphe dont les sommets sont dans $\chi^{-1}([0,+\infty[)$, et dont les ar\^etes sont toutes les ar\^etes joignant de tels sommets n'est pas connexe, ceci  pour un syst\`eme g\'en\'erateur fini quelconque de $G$. Cette non-connexit\'e est ind\'ependante de ce syst\`eme g\'en\'erateur). L'existence d'un caract\`ere ``exceptionnel" lorsque $DG$ n'est pas de type fini est [B-S, theorem 4.2], l'un des r\'esultats de base de la th\'eorie.

Soit $d:Z\to X$ une d\'esingularisation, qui est un isomorphisme au-dessus de $(X-Supp(\Delta))$. Elle induit, par [Ko 93,7.2-7.9], un isomorphisme au niveau des groupes fondamentaux.
On d\'eduit de $\chi$ une fonction $h:Z_{ab}\to \bR$ sur le rev\^etement ab\'elien (sans torsion) maximal de $Z$, de groupe $\pi_1(X)_{ab}/Torsion$. La d\'eriv\' ee de cette fonction descend en une $1$-forme r\'eelle $d$-ferm\'ee $w$, partie r\'eelle d'une $1$-forme holomorphe $\omega$ sur $Z$. 

 Soit alors $r:\widetilde{(X\vert \Delta)}\to X$ le rev\^etement universel de $(X\vert \Delta)$, et $\tilde d:\widetilde Y\to \widetilde{(X\vert \Delta)}$ une d\'esingularisation $G$-\' equivariante qui domine $d:Z\to X$, et qui est un isomorphisme au-dessus de $(X-Supp(\Delta))$. La fonction $h$ se rel\`eve en une fonction $\tilde h:\widetilde Y\to \bR$ telle que $\widetilde{h}^{-1}([s,+\infty[)$ ne soit pas connexe pour $s>0$ assez grand (par le lemme 2.2 de [De 06], adapt\'e). Plus pr\'ecis\'ement, on peut exprimer topologiquement $(X-Supp(\Delta))$ comme un $CW$-complexe dont le $1$-squelette est un bouquet de cercles (ne rencontrant donc pas $Supp(\Delta))$. Son image r\'eciproque dans $\widetilde{(X\vert \Delta)}$ est alors topologiquement le graphe de Cayley de $G$ (pour le syst\`eme non minimal de g\'en\'erateurs d\'efini par ce bouquet de cercles).

Alors $\widetilde{h}^{-1}([s,+\infty[)$ n'est pas connexe, pour $s>0$ assez grand, puisque si $D\subset \widetilde{(X\vert \Delta)}$ est un domaine fondamental compact pour $G$, $\widetilde h$ est born\'ee sur $D$, et que $\widetilde h(g.D)=log (\chi(g))+\widetilde h(D)$.

Soit alors  $\widetilde w$ le rel\`evement \`a $\widetilde Y$ de $w$. Le second argument (section $2.2.2$ de [De 06], qui en est le point crucial) s'applique encore directement \`a notre situation, et montre que le feuilletage de $\widetilde{Y'}$ d\'efini par la $1$-forme holomorphe $g^*(\omega)=d\widetilde{F}$ a une feuille $\widetilde{F}^{-1}(\psi)$ dont une composante connexe $\widetilde{C_{\omega}}$ est critique, c'est-\`a-dire telle que $d\widetilde{F}=0$ sur cette composante (de codimension un). L'image dans $X$ de cette composante est une feuille compacte du feuilletage de $X$ d\'efini par $\omega$ puisque les fibres de $\widetilde{d}$ sont envoy\'ees sur des sous-ensembles analytiques de codimension complexe au moins $2$, et que $r$ est finie au voisinage de chacun des points de $\widetilde{(X\vert \Delta)}$. (Remarquer que cette feuille compacte ne provient \'evidemment pas, en g\'en\'eral, d'une feuille critique de l'image r\'eciproque de $\omega$ sur le rev\^etement universel de $X$. L'annulation additionnelle provient de la ramification de $r$ le long des composantes de $\Delta)$.

Il r\'esulte alors de [DG 05, 4.1] ou de [Si 93-2] que les feuilles de $\omega$ sont toutes compactes et d\'efinissent une fibration $f:X\to C$ sur une courbe de genre $g\geq 1$ telle que $\omega$ provienne de $C$. La base orbifolde $(C\vert\Delta_{f,\Delta})$ est alors n\'ecessairement hyperbolique, puisque l'on a une suite exacte de groupes ([Ca 07, Prop.11.7]): $$1\to \pi_1(F\vert \Delta_F)\to G\to \pi_1(C\vert \Delta_{f,\Delta})\to 1,$$ dans laquelle $(F\vert \Delta_F)$ est la fibre orbifolde (g\'en\'erique) de $f: (X\vert \Delta)\to C$, de sorte que $\pi_1(F\vert \Delta_F)$ est de type fini, la classe exceptionnelle $\chi$ provenant de $C$. En effet si $(C\vert \Delta_C)$ n'est pas hyperbolique, $C$ est elliptique et $\Delta_C=0$.  On a donc une suite exacte (d\'eduite de la pr\'ec\'edente): $$1\to \pi_1(F\vert \Delta_F)\to Ker (\chi\circ f_*)\to Ker(\chi)\to 1$$

Mais $Ker(\chi)\subset \bZ^{\oplus 2}$ est de type fini, donc aussi $Ker(\chi\circ f_*)$, puisque $\pi_1(F\vert \Delta_F)$ est de type fini. Ceci contredit le fait qu'il existe un $\chi$ exceptionnel provenant de $C$ ([B-S, theorem 4.2]). $\square$

\begin{corollary}\label{qr} Soit $(X\vert \Delta)$ une orbifolde lisse et K\" ahler dont aucun rev\^etement orbifolde-\'etale n'admet de morphisme orbifolde sur une courbe orbifolde  hyperbolique (condition satisfaite si, par exemple, $\pi_1(X\vert \Delta)$ est r\'esoluble), alors tout quotient r\'esoluble de $\pi_1(X\vert \Delta)$ est virtuellement nilpotent.
\end{corollary}

\begin{re}\label{descente} Si $(X\vert \Delta)$ une orbifolde lisse et K\" ahler, et $u:(X'\vert \Delta')\to (X\vert \Delta)$ un rev\^etement orbifolde \'etale fini qui admet un morphisme orbifolde $f':(X'\vert \Delta')\to (C'\vert \Delta_{C'})$ sur une courbe orbifolde hyperbolique. Alors $(X\vert \Delta)$ n'admet pas toujours de morphisme orbifolde sur une courbe orbifolde hyperbolique. On ne peut donc simplifier les \'enonc\'es \ref{resnilp} et \ref{dgti} en supprimant les conditions sur les rev\^etements. Consid\'erons par exemple $X_0$ la jacobienne d'une courbe g\'en\'erique $C$ de genre $2$, $X_0$ identifi\'ee aux classes $[D]$ de diviseurs $D$ de degr\'e $2$ modulo \'equivalence lin\'eaire sur $C$. Soit $g:X\to X_0$ l'\'eclat\'e de $X_0$ au-dessus du point $[K_C]$. On a un morphisme naturel $u:X':=C\times C\to X$, fini et de degr\'e $2$, ramifi\'e le long de la diagonale de $X'$, et au-dessus de la transform\'ee stricte $C_1$ de $C\subset X_0$ (par plongement Jacobien par points doubles $2.y,y\in C)$. On ne met pas de structure orbifolde sur $X'$, mais on munit $X$ du diviseur orbifolde $\Delta:=(1-\frac{1}{2}). C_1$. Le morphisme $u$ est donc orbifolde-\'etale. Il est clair que $X'$ a des fibrations sur des courbes hyperboliques. Cependant, $(X\vert (1-\frac{1}{2}). C_1)$ n'a pas de morphisme sur une courbe orbifolde hyperbolique, puisque s'il en existait on aurait (par [Ca 07, th\'eor\`eme 8.17]) un morphisme presque-holomorphe, donc holomorphe $\varphi: X\to B$ sur une courbe. Puisque $X_0$ n'a pas de tel morphisme (sinon $X_0$ serait isog\`ene \`a un produit de courbes elliptiques, contredisant la g\'en\'ericit\'e de $C)$, $\varphi$ doit avoir sur $X_0$ un point d'ind\'etermination simple unique au point $[K_C]$, et les diviseurs fibres un nombre d'auto-intersection \'egal \` a $1$, ce qui contredit la parit\'e de la forme d'intersection sur $X_0$.
\end{re}

\section{Morphisme d'Albanese.}\label{alb}

\subsection{Factorisation de Stein du morphisme d'Albanese.}\label{albstein}

Soit $X$ compacte K\" ahler et connexe. Soit $\alpha_X:X\to Alb(X)$ son morphisme d'Albanese, $Z$ un mod\`ele lisse de $\alpha_X(A)\subset Alb(X)$, et enfin $Y$ un mod\`ele lisse de la factorisation de Stein de $\alpha_X:X\to Z$ (que l'on supposera holomorphe, ainsi que $f:X\to Y$, quitte \`a modifier $X$, ce qui ne change pas $\pi_1(X))$. 

On obtient ainsi deux applications holomorphes: $f:X\to Y$ et $g:Y\to Z$, avec $g\circ f:X\to Z$ bim\'eromorphe \`a $\alpha_X: X\to \alpha_X(X)$, avec $f$ \`a fibres connexes et $g$ g\'en\'eriquement finie, telles que $Alb(Z)\cong Alb(X)$.

\

Soit maintenant $(X\vert \Delta)$ une structure orbifolde lisse sur $X$. Elle induit sur $Y$ et $Z$ des structures orbifoldes $\Delta_Y:=\Delta_{f,\Delta}$ et $\Delta_Z:=\Delta_{g\circ f,\Delta}$. On supposera (par modification ad\'equate de $X,Y,Z)$ que les morphismes $f$ et $g\circ f$ sont {\it nets} au sens de [Ca 07]. Ceci implique ([Ca 07]) que pour tout morphisme $k:Z\to B$ sur une courbe projective lisse $B$, si $\Delta_C:=\Delta_{k,\Delta_Z}$, alors: $\Delta_{k\circ g\circ f,\Delta}:=\Delta_{k,\Delta_Z}$ (et, de mani\`ere similaire: $\Delta_{k\circ g\circ f,\Delta}:=\Delta_{k\circ g,\Delta_Y})$. Remarquons que tout morphisme de $X$ sur une courbe de genre $g>0$ se factorise par $Z$, par universalit\'e de $Alb(X)$.

\

Soit alors $f_j:(X\vert \Delta)\to (C_j\vert \Delta_j)$ l'ensemble (fini, par le th\'eor\`eme de Kobayashi-Ochiai, voir [A 92] pour ce cas particulier) des morphismes orbifoldes surjectifs et \`a fibres connexes sur des orbifoldes de courbes hyperboliques de genre $g>0$, notant pour tout $j$:$\Delta_j:=\Delta_{f_j,\Delta}$. 

On notera enfin: $h^X:X\to \Pi_jC_j$ le morphisme produit des $f_j$, et $V\subset \Pi_j C_j$ son image d\'esingularis\'ee. On notera encore $h^X: X\to V$ le morphisme r\'esultant (apr\`es modification de $X)$. On supposera encore que $h: X\to V$ est {\it net}. On munira enfin $V$ de la structure orbifolde $\Delta_V:=\Delta_{h,\Delta}$. Par la propri\'et\'e de factorisation rappel\'ee ci-dessus, il existe donc un unique morphisme $h:Z\to V$ tel que $h^X=h\circ g\circ f$.

Nous obtenons ainsi une compos\'ee:

\centerline{
\xymatrix{ (X\vert \Delta)\ar[r]^{f}&(Y\vert \Delta_Y)\ar[r]^{g}&(Z\vert \Delta_Z)\ar[r]^{h}&(V\vert \Delta_V)\\
}}

\subsection{Green-Lazarsfeld  et morphisme d'Albanese.} \label{glalb}

\

Si $(X\vert \Delta)$ est une orbifolde K\" ahl\'erienne, on note $\Sigma_c^1(X\vert \Delta)$ la partie continue de $\Sigma^1(X\vert \Delta)$, de sorte que $\Sigma^1(X\vert \Delta)$ est r\'eunion de $\Sigma_c^1(X\vert \Delta)$ et d'un nombre fini de points isol\'es (qui sont de torsion).

\begin{proposition}\label{glalb'} Dans la situation et avec les notations du \S\ref{albstein} pr\'ec\'edentes, les applications naturelles d\'eduites de $f,g,h$ induisent des bijections: 

$\Sigma_c^1(V\vert \Delta_{V})\to\Sigma_c^1(Z\vert \Delta_{Z})\to\Sigma_c^1(Y\vert \Delta_{Y})\to\Sigma_c^1(X\vert \Delta_{X})$.

\end{proposition}

{\bf D\'emonstration:} Il suffit, par fonctorialit\'e de l'\'etablir pour la compos\'ee. Dans ce cas l'assertion r\'esulte ce ce que $\Sigma_c^1(X\vert \Delta_{X})$ est la r\'eunion des $f_j^*(\Sigma^1(C_j\vert \Delta_{j}))$, du fait que  les $f_j$ se factorisent par $g\circ f$, et de l'\'egalit\'e: $\Delta_{h\circ g\circ f,\Delta}:=\Delta_V$ $\square$

\

Ces \'egalit\'es sugg\`erent de comparer les quotients r\'esolubles des groupes fondamentaux des $4$ orbifoldes pr\'ec\'edentes.

\subsection{Quotients r\'esolubles  et morphisme d'Albanese.} \label{qralb}

\

Les applications $f:X\to Y, $ $g:Y\to Z$ et $h: Z\to V$ induisent des morphismes de groupes naturels  $f_*: G:=\pi_1(X\vert \Delta)\to \pi_1(Y\vert \Delta_Y):=G_Y$,  $g_*: G_Y:=\pi_1(Y\vert \Delta_Y)\to \pi_1(Z\vert \Delta_Z):=G_Z$, et  $h_*: G_Z=\pi_1(Z\vert \Delta_Z)\to \pi_1(V\vert \Delta_V):=G_V$. 

\

On notera aussi: $G_f=Ker(f_*)$, $G_{gf}:=Ker((g\circ f)_*)$, et $G_{hgf}:=Ker((h\circ g\circ f)_*)$. Ce sont des sous-groupes normaux de $G$. 

On a bien s\^ur des inclusions: $G_f<G_{gf}<G_{hgf}$. 

Remarquer que $G_{hgf}=G$ lorsque $\Sigma_c^1(X\vert \Delta)=\emptyset$ ou, de mani\`ere \'equivalente, lorsque $(X\vert \Delta)$ n'a pas de fibration sur une courbe orbifolde hyperbolique de genre $g>0$.

\

Un morphisme de groupes surjectif $\varphi:G\to R$, avec $R$ r\'esoluble sera appel\'e {\it un quotient r\'esoluble de $G$}. Les images $\varphi(G_f)$, $\varphi(G_{gf})$, et $\varphi(G_{hgf})$ sont donc alors des sous-groupes (r\'esolubles) normaux de $R$.

\begin{question}\label{qqralb} Avec les notations de \ref{albstein} ci-dessus, soit $\varphi:G\to R$ un morphisme de groupes surjectif, avec $R$ r\'esoluble (on dira aussi que $R$ est un quotient r\'esoluble de $G)$.

1. Les groupes $\varphi(G_f)$ et $\varphi(G_{gf})$ sont-ils finis?

2. Le groupe $\varphi(G_{hgf})$ est-il virtuellement nilpotent?

3. Tout quotient r\'esoluble de $(G_{hgf})$ est-il virtuellement nilpotent?

\end{question}

Une r\'eponse affirmative \`a la question \ref{qqralb}.1 signifierait donc que les quotients r\'esolubles de $G$ et de $G_Z$ sont les m\^emes, \`a commensurabilit\'e pr\`es. Et donc que ces quotients sont ceux des (groupes fondamentaux) orbifoldes de sous-vari\'et\'es de tores complexes. Il est d\'emontr\'e dans [Ca 95] que c'est le cas pour les quotients nilpotents de $G$.

Nous r\'epondrons (affirmativement) \`a ces questions dans le cas de $G_f$ seulement.

\begin{theorem}\label{rqralb} Avec les notations pr\'ec\'edentes, soit $(X\vert \Delta)$ une orbifolde lisse et K\" ahler, et $\varphi:G\to R$, un quotient r\'esoluble de $G:=\pi_1(X\vert \Delta)$. Alors:

1. $\varphi(G_f)$ est fini.

2. Si $(X\vert \Delta)$ n'a pas de fibration sur une courbe hyperbolique, le groupe $\varphi(G_{hgf})=R$ est virtuellement nilpotent. 
\end{theorem}

\begin{remark} 
La finitude de $\varphi(G_{gf})$ est plus d\'elicate que celle de $\varphi(G_f)$, car non pr\'eserv\'ee {\it a priori} par rev\^etements orbifolde-\'etales.
\end{remark}

{\bf D\'emonstration:} Assertion 2. Si $(X\vert \Delta)$ n'a pas de fibration sur une courbe hyperbolique, alors $G_{hgf}=G$, et tout quotient r\'esoluble de $G$ est virtuellement nilpotent par \ref{qr}, ce qui \'etablit la seconde partie de la seconde assertion. La finitude de $\varphi(G_f)$ dans ce cas est un cas particulier de l'assertion 1, d\'emontr\'ee ci-dessous, puisque $R=\varphi(G_{hgf})$ est alors virtuellement nilpotent, donc r\'esiduellement fini.

\

La d\'emonstration de l'assertion 1, que nous abordons maintenant, est nettement plus d\'elicate. Elle utilise les r\'esultats et d\'efinitions interm\'ediaires \ref{qralb'}-\ref{sectexist} qui suivent:

\begin{theorem}\label{qralb'} (Voir [Ca 01, th\'eor\`eme 4.1])\footnote{Dans [Ca 01, 4.1], on affirme que $\varphi(G_f\cap G')$ est lui-m\^eme fini. La d\'emonstration \'etablit d'abord la finitude de l'ab\'elianis\'e, et en d\'eduit la finitude de ce groupe par l'interm\'ediaire du lemme [Ko 93: J. Koll\`ar. Inv. Math.113 (1993), 177-215, Prop. 6.4] dont la d\'emonstration est cependant incompl\`ete. Nous n'\' enon\c cons donc dans \ref{qralb'} ci-dessus que la partie de ce r\'esultat qui ne d\'epend pas de [Ko 93, 6.4], et \'etablissons donc [Ca 01, 4.1] int\'egralement dans \ref{qralb} ci-dessus.}
 Soit $f:(X\vert \Delta)\to (Y\vert \Delta_Y)$ l'application d\'efinie ci-dessus, et $\varphi:G\to R$ un quotient r\'esoluble de $G$. 
 
 Alors, pour tout sous-groupe $G'$ d'indice fini de $G$: $(\varphi(G_f\cap G'))_{ab}$ est un groupe fini.
\end{theorem}

{\bf D\'emonstration:}  Soit $(F\vert \Delta_F)$ l'une des fibres orbifoldes g\'en\'eriques de $f: (X\vert \Delta)\to (X\vert \Delta)$. Alors $G_f\subset G$ est aussi l'image de $\pi_1(F\vert \Delta_F)$ dans $G:=\pi_1(X\vert \Delta)$ (par [Ca 07, 11.2]). 
Le r\'esultat est donc \'etabli lorsque $\Delta=\Delta_Y=0$ dans [Ca 01, 4.1]. La d\'emonstration de loc. cit. s'adapte en fait directement au cas plus g\'en\'eral consid\'er\'e ici en remarquant que $(\pi_1(F\vert \Delta_F))_{ab}/Torsion=(\pi_1(F))_{ab}/Torsion$ $\square$

\

Le probl\`eme crucial de la d\'emonstration de \ref{qralb} r\'eside dans la non finitude r\'esiduelle des groupes de type fini r\'esolubles g\'en\'eraux. Nous aurons besoin d'une version relative.

\begin{definition}\label{resfin} Soit $L<G$ un sous-groupe normal. Nous dirons que $L$ est {\it r\'esiduellement fini dans $G$} si, pour tout sous-groupe $L'<L$ d'indice fini, il existe un sous-groupe $G'<G$ d'indice fini tel que: $(G\cap L)\subset L'$.

\end{definition}

\begin{example}\label{exresfin}

\

1. Si $G$ est r\'esiduellement fini, tout $L<G$ est r\'esiduellement fini dans $G$.

2. Si $G=(G/L)\ltimes L$ est produit semi-direct de $L$, $L$ de type fini, et d'un sous-groupe $(G/L)<G$, alors $L$ est  r\'esiduellement fini dans $G$. 

En effet: si $L'$ est caract\'eristique dans $L$, alors $G':=(G/L)\ltimes L'$ r\'esoud le probl\`eme pour $L'$. Or tout sous-groupe d'indice fini de $L$ contient un sous-groupe caract\'eristique de $L$ d'indice fini, puisque $L$ est de type fini.

\end{example}

\begin{lemma}\label{qra}  Dans la situation de \ref{rqralb}, $\varphi(G_f)$ est fini si $G_f$ est r\'esiduellement fini dans $G$.\end{lemma}

{\bf D\'emonstration (de \ref{qra}):} Elle d\'ecoule directement de \ref{qralb'} ci-dessus, puisque, d'apr\`es [Ca 01, 3.6.1] appliqu\'e \`a $S=\varphi(G_f)$, un groupe r\'esoluble $S$ est fini si $S'_{ab}$ est fini pour tout sous-groupe $S'<S$ d'indice fini $\square$

\

Nous allons maintenant nous ramener au cas o\`u $G_f$ est r\'esiduellement fini dans $G$.

\begin{lemma}\label{sect} Soit $f:(X\vert \Delta)\to (Y\vert \Delta_Y)$, avec $\Delta_Y:=\Delta_{f,\Delta}$, $(F\vert \Delta_F)$ d\'esignant la fibre orbifolde g\'en\'erique de $f$, et $f$ {\it nette}. Si $f$ admet une section orbifolde divisible $\sigma:(Y\vert \Delta_Y)\to (X\vert \Delta)$, alors la suite exacte naturelle:$$\pi_1(F\vert \Delta_F)\to \pi_1(X\vert \Delta)\to \pi_1(Y\vert \Delta_Y)\to 1$$ est scind\'ee, et l'image $L=G_f$ de $\pi_1(F\vert \Delta_F)$ dans $G$ est r\'esiduellement finie dans $G:=\pi_1(X\vert \Delta_X)$.

Rappelons [Ca 07, \S5] que $\sigma$ est une section orbifolde divisible de $f$ si c'en est une section (au sens usuel), et un morphisme orbifolde divisible.
\end{lemma}

{\bf D\'emonstration:} Le morphisme orbifolde $\sigma$ induit [Ca 07, \S11] un morphisme de groupes $\sigma_*:\pi_1(Y\vert \Delta_Y)\to \pi_1(X\vert \Delta)$ dont l'image $\Sigma$ est un suppl\'ementaire de $G_f$ dans $\pi_1(X\vert \Delta)$, qui est donc produit semi-direct de $\Sigma$ et de $G_f$. Donc $G_f$ est r\'esiduellement fini dans $G$ $\square$

\begin{lemma}\label{sectexist} Soit $f:(X\vert \Delta)\to (Y\vert \Delta_Y)$ comme ci-dessus. Il existe alors un diagramme commutatif de morphismes orbifoldes (divisibles): 

\

\
\centerline{
\xymatrix{(X'\vert \Delta')\ar[r]^{u}\ar[d]_{f'} &(X\vert \Delta)\ar[d]_{f}\\(Y'\vert \Delta_{Y'})\ar[r]^v&(Y\vert \Delta_Y)\\
}}

\

tel que:

1. Les vari\'et\'es $X',Y'$ sont K\" ahler.

2. Les fibres orbifoldes g\'en\'eriques de $f$ et $f'$ coincident (ie: $f'$ est g\'en\'eriquement d\'eduite de $f$ par le changement de base $v)$.

3. $f':(X'\vert \Delta')\to (Y'\vert \Delta_{Y'})$ admet une section orbifolde $\sigma$.

4. Dans cette situation, $Alb(X')=Alb(X)\times_{Alb(Y)}Alb(Y')$, et les images de $F'$ dans $Alb(X')$ et de $F$ dans $Alb(X)$ ont donc la m\^eme dimension.

\end{lemma}

{\bf D\'emonstration:} On choisit: $Y':=X$. Il existe donc une section tautologique (diagonale) naturelle $\sigma: Y'\to X':=Y'\times _Y X$ de la fibration $f': X'\to Y'$. Nous noterons encore $X'$ un mod\`ele lisse ad\'equat de $Y'\times _Y X$. On choisit alors sur $X'$ une structure orbifolde g\'eom\'etrique $\Delta'$ telle que:

1. La fibre orbifolde g\'en\'erique $(F'\vert\Delta_{F'})$ soit \'egale \`a son image $(F\vert \Delta_F)$. 

2. Le morphisme naturel $u:(X'\vert \Delta')\to (X\vert \Delta)$ soit un morphisme orbifolde: c'est toujours possible en augmentant (au sens divisible) les multiplicit\'es sur les composantes $f'$-verticales du support de $\Delta'$.

3. La fibration $f'$ est nette.

On d\'efinit alors $\Delta_Y':=\Delta_{f',\Delta'}$. 

En g\'en\' eral, la section $\sigma:(Y'\vert \Delta_{Y'})\to (X'\vert \Delta')$ n'est pas un morphisme orbifolde  (la $\Delta_{Y'}$-multiplicit\'e d'un diviseur $E'$ irr\'eductible effectif de $Y'$ peut diviser strictement celle des composantes $D'$ de $\Delta'$ que rencontre $\sigma$ au-dessus de $E')$. On peut n\'eammoins, pour tout tel $E'$, augmenter au sens divisible les $\Delta'$-multiplicit\'es des composantes irr\'eductibles de $(f')^*(E')$ de telle sorte que la $\Delta_{Y'}$-multiplicit\'e de $E'$ soit \'egale \`a la plus grande des $\Delta'$-multiplicit\'es des composantes $D'$. La section $\sigma:(Y'\vert \Delta_{Y'})\to (X'\vert \Delta')$ est alors un morphisme orbifolde divisible.

4. Ceci r\'esulte de la commutativit\'e du diagramme:

\
\centerline{
\xymatrix{H^1(X)\ar[r]^{u^*}\ar[d]_{f^*} &H^1(X')\ar[d]^{(f')^*}\\H^1(F)\ar[r]^{u^*}&H^1(F')\\
}}

\

dans lequel les coefficients sont rationnels $\square$

\

Nous pouvons maintenant d\'emontrer l'assertion (1) du th\'eor\`eme \ref{rqralb}:

\

{\bf D\'emonstration de \ref{rqralb}.(1):} Elle r\'esulte imm\'ediatemment de \ref{qra} et des lemmes \ref{sect} et \ref{sectexist} ci-dessus: en effet (avec les notations de \ref{sectexist}), \ref{qra} et \ref{sect} montrent que l'assertion est vraie pour $(X'\vert \Delta')$, et que l'assertion pour $(X'\vert \Delta')$ l'implique pour $(X\vert \Delta)$ aussi. $\square$

\section{Orbifoldes sp\'eciales et ``Coeur".}\label{orbspec}

\subsection{Rappels sur les orbifoldes sp\'eciales.}\label{ros}

Nous allons rappeler, suivant [Ca 04] et [Ca 07], leur place centrale dans la classification. Il n'est pas possible de donner ici, faute de place, les d\'efinitions n\'ecessaires. Nous n'\'evoquons donc ici que la structure des constructions et r\'esultats, en laissant de c\^ot\'e certaines nuances et pr\'ecisions. 

\

 Pour simplifier la description, nous consid\'erons tout d'abord le cas d'une vari\'et\'e (K\"ahler compacte) $X$ sans structure orbifolde. Les r\'esultats et conjectures (mais non les d\'efinitions) s'\'etendent sans changement aux $(X\vert \Delta)$ ``lisses".

\

Rappelons que si $X$ est une vari\'et\'e K\" ahl\'erienne compacte et connexe, on a introduit dans [Ca 04] la notion de {\it base orbifolde} $(Y\vert \Delta_f)$ d'une fibration m\'eromorphe dominante $f:X\dasharrow Y$: il s'agit d'un rev\^etement ramifi\'e virtuel de $Y$ \'eliminant en codimension $1$ les fibres multiples de $f$, et constitu\'e du couple $Y$ et d'un diviseur orbifolde $\Delta_f$, diviseur de ramification du rev\^etement ramifi\'e virtuel associ\'e.

\

Une nouvelle classe de vari\'et\'es, dites {\it sp\'eciales}, est ensuite introduite. Ces vari\'et\'es sont d\'efinies par les deux conditions \'equivalentes A et B suivantes: 

{\bf A.} Pour tout $p>0$ et tout sous-faisceau $L\subset \Omega^p_X$, coh\'erent de rang $1$, on a: $\kappa(X,L)<p$. (En particulier, $X$ n'a pas de fibration m\'eromorphe $f$ sur une vari\'et\'e $Y$ de type g\'en\'eral et de dimension $p>0$, sans quoi $f^*(K_Y)$ fournirait une contrediction).

{\bf B.} Aucune fibration m\'eromorphe dominante $f:X\dasharrow Y$ n'a de base orbifolde $(Y\vert\Delta_{f})$ de type g\'en\'eral.

Les vari\'et\'es rationnellement connexes, et les vari\'et\'es avec $\kappa=0$, sont sp\'eciales. Cette classe n'est cependant pas caract\'eris\'ee par la dimension de Kodaira: les surfaces avec $\kappa=1$ peuvent \^etre soit sp\'eciales, soit non-sp\'eciales. Plus pr\'ecis\'ement, une surface est sp\'eciale si et seulement si $\kappa\leq 1$ et $\pi_1$ est presque-ab\'elien, ce qui montre que cette propri\'et\'e est stable par d\'eformation. En dimension $3$ ou plus, la situation est beaucoup plus compliqu\'ee. Conditionnellement en une version orbifolde de la conjecture $C_{n,m}$, on montre dans [Ca 07,\S 10] que les vari\'et\'es (et orbifoldes) sp\'eciales se d\'ecomposent canoniquement en tours de fibrations \`a fibres orbifoldes g\'en\'erales ayant soit $\kappa=0$, soit $\kappa_+=-\infty$ (une version faible de la connexit\'e rationnelle). Voir \ref{crj} pour plus de d\'etails.

Cette d\'ecomposition am\`ene naturellement \`a  conjecturer (en particulier) que:

{\bf C.} Si $X$ est sp\'eciale, $\pi_1(X)$ est presque-ab\'elien.

{\bf D.} $X$ est sp\'eciale si et seulement si sa pseudo-m\'etrique de Kobayashi est nulle, si et seulement si elle est potentiellement dense (lorsque d\'efinie sur un corps de nombres, disons).




 Apr\`es ce rappel, nous revenons au cas g\'en\'eral d'une orbifolde lisse $(X\vert \Delta)$. On s'attend \`a ce que les propri\'et\'es soient analogues au cas o\`u $\Delta=0$, et l'objectif du pr\'esent texte est justement de le montrer pour certaines propri\'et\'es du groupe fondamental. Nous renvoyons \`a [Ca 04] et [Ca 07] pour plus de d\'etails sur les notions utilis\'ees ici.

Soit $(X\vert\Delta)$ une orbifolde lisse et K\" ahler. Rappelons ([Ca 07, 8.1]) que $(X\vert\Delta)$ est dite {\it sp\'eciale} s'il n'existe pas de fibration m\'eromorphe (au sens orbifolde) dominante $g:(X\vert\Delta)\dasharrow (Y\vert\Delta_Y)$ sur une orbifolde $(Y\vert\Delta_Y)$ de type g\'en\'eral avec $dim(Y)=\kappa(Y\vert\Delta_Y)>0$. Une d\'efinition \'equivalente est la version orbifolde de la d\'efinition {\bf A.} donn\'ee ci-dessus lorsque $\Delta=0$.

 Les exemples fondamentaux de telles orbifoldes sont les $(X\vert\Delta)$ qui sont soit rationnellement connexes, soit avec $\kappa(X\vert\Delta)=0$. Une fibration dont la fibre orbifolde g\'en\'erique et la base orbifolde sont sp\'eciales est sp\'eciale (une propri\'et\'e fausse dans la cat\'egorie des vari\'et\'es sans structure orbifolde). Les tours orbifoldes de telles fibrations sont donc sp\'eciales. La r\'eciproques \'etant vraie conditionnellement en $C_{n,m}^{orb}$ (Voir \ref{crj}). Si $a_X:X\to W$ est la r\'eduction alg\'ebrique de $X$, et $(X\vert \Delta)$ une structure orbifolde sur $X$, la fibre orbifolde g\'en\'erique $(F\vert \Delta_F)$ de $a_X$ est sp\'eciale (ceci est d\'emontr\'e dans [Ca 04] lorsque $\Delta=0$, mais la d\'emonstration s'adapte sans changement au cas g\'en\'eral). En particulier: $(X\vert \Delta)$ est sp\'eciale si $a(X)=0$, $a(X)$ \'etant la dimension alg\'ebrique de $X$.

 \subsection{Orbifoldes sp\'eciales et morphisme d'Albanese.}

 Soit  $(X\vert\Delta)$ une orbifolde lisse et K\" ahler, et $\alpha_X:X\to Alb(X)$ son morphisme d'Albanese.

 \begin{theorem}\label{specalb} Si $(X\vert\Delta)$ est sp\'eciale, alors: 
 
 0. $X$ est sp\'eciale.

1. $\alpha_{X}:X\to Alb(X)$ est surjective et connexe.

2. $\Delta_{(\alpha_X,\Delta)}=0$.

3. Tout rev\^etement orbifolde-\'etale  fini $(X'\vert\Delta')\to(X\vert\Delta)$ est sp\'ecial.

 \end{theorem}

{\bf D\'emonstration:} 0. Il existe en effet un morphisme orbifolde naturel surjectif: $(X\vert \Delta)\to X=(X\vert 0)$, et l'image d'une orbifolde sp\'eciale l'est aussi. L'assertion 1 est [Ca 04, prop. 5.3], puisque $X$ est sp\'eciale. L'assertion 2 se d\'emontre exactement comme l'\'enonc\'e [Ca 04, prop. 5.3], qui affirme que si $X$ est sp\'eciale, alors $\alpha_X$ est surjective, connexe et sans fibre multiple en codimension $1$. (La d\'emonstration s'applique sans changement lorsque $\Delta\neq 0)$.  

L'assertion 3 r\'esulte de [Ca 07, prop. 9.11] $\square$

\subsection{Quotients r\'esolubles du $\pi_1$ des orbifoldes sp\'eciales.}

La conjecture centrale concernant le groupe fondamental de ces orbifoldes g\'eom\'etriques sp\' eciales est (voir [Ca 07]):

\begin{conjecture}\label{cj} Si $(X\vert\Delta)$ est sp\'eciale, alors $\pi_1(X\vert\Delta)$ est virtuellement ab\'elien.
\end{conjecture}

Nous allons v\'erifier cette conjecture pour les quotients r\'esolubles de $\pi_1(X\vert\Delta)$, et pour une classe un peu plus large d'orbifoldes.

\begin{theorem}\label{albsurj} Soit $(X\vert\Delta)$ une orbifolde lisse et K\" ahler. Si, pour tout rev\^etement orbifolde-\'etale fini $r: (X'\vert\Delta')\to (X\vert\Delta)$, l'application d'Albanese $\alpha_{X'}:X'\to Alb(X')$ est surjective\footnote{Une condition plus faible (\'egalement n\'ecessaire) est en fait suffisante: l'injectivit\'e de l'application $(\alpha_{X'})^*:H^2(Alb(X'),\bQ)\to H^2(X',\bQ)$. Voir [Ca 95]: la d\'emonstration est la m\^eme.}, alors tout quotient virtuellement r\'esoluble de $\pi_1(X\vert\Delta)$ est virtuellement ab\'elien.
\end{theorem}

\begin{remark} Si l' orbifolde $(X\vert\Delta)$, lisse et K\" ahler, est sp\'eciale, alors ses rev\^etements orbifolde-\'etales finis le sont aussi, d'apr\`es \ref{specalb}.(3).
Le th\'eor\`eme \ref{albsurj} s'applique donc, d'apr\`es \ref{specalb}.(1), en particulier, \`a toute telle orbifolde sp\'eciale.
\end{remark}

{\bf D\'emonstration:} Hypoth\`ese et conclusion sont pr\'eserv\'ees par rev\^etements orbifolde-\'etale. L'hypoth\`ese implique \'evidemment que $(X\vert\Delta)$ n'a pas de fibration sur une courbe orbifolde hyperbolique. Les quotients r\'esolubles de $\pi_1(X\vert\Delta)$ sont donc virtuellement nilpotents. Il suffit donc (apr\`es rev\^etement orbifolde-\'etale ad\'equat) de montrer que si un tel quotient $R$ est nilpotent et sans torsion, il est virtuellement ab\'elien. Mais alors $R$ est aussi un quotient de $\pi_1(X)$, et le r\'esultat principal de [Ca 95] implique l'assertion. $\square$

\begin{corollary} Soit $(X\vert\Delta)$ une orbifolde lisse et K\" ahler sp\'eciale, et $R$ un groupe quotient de $\pi_1(X\vert\Delta)$. Alors $R$ est virtuellement ab\'elien dans les deux cas suivants:

1. $R$ est virtuellement r\'esoluble.

2. $R<Gl(n,\bC)$ est lin\'eaire.

De plus, si $q(X'):=h^0(X',\Omega^1_{X'})=0$ pour tout rev\^etement \'etale fini $X'$ de $X$, alors $R$ est fini.

\end{corollary}

{\bf D\'emonstration:} Si $R$ est virtuellement r\'esoluble, il suffit, par \ref{albsurj} de v\'erifier que l'application d'Albanese de $X$ est surjective, puisque les rev\^etements orbifolde-\'etales finis de $(X\vert\Delta)$ sont sp\'eciaux. Mais $X$ \'etant elle-m\^eme sp\'eciale, cette surjectivit\'e r\'esulte de [Ca 04]. La premi\`ere assertion est ainsi \'etablie. 

Si $R$ est lin\'eaire, il admet, par le lemme de Selberg, un sous-groupe $R'$ d'indice fini sans torsion, qui est donc un quotient de $\pi_1(X)$ (puisque le noyau de $\pi_1(X\vert\Delta)\to \pi_1(X)$ est engendr\'e par des \'el\'ements de torsion). Puisque $X$ est aussi sp\'eciale, l'assertion r\'esulte encore de [Ca 04, 7.7].

Pour d\'emontrer la derni\`ere assertion, supposons (apr\`es rev\^etement orbifolde-\'etale fini ad\'equat) $R$ ab\'elien infini sans torsion. Alors $R$ est un quotient de $(\pi_1(X\vert\Delta))_{ab}$, et donc aussi de $(\pi_1(X))_{ab}$, puisque le noyau de  $(\pi_1(X\vert\Delta))\to (\pi_1(X))$ est engendr\'e par des \'el\'ements de torsion. Contradiction, puisque $q(X)=0$, par hypoth\`ese. $\square$



\subsection{Le ``Coeur".}\label{redalge}

Rappelons que l'on montre dans [Ca 07] (et [Ca 04] lorsque $\Delta=0)$, l'existence et l'unicit\'e, pour toute orbifolde lisse  $(X\vert \Delta)$, avec $X$ compacte et connexe K\" ahler, d'une nouvelle fibration $c:(X\vert \Delta)\to W$, presque-holomorphe, (appel\'ee ``le coeur" de $(X\vert \Delta)$) telle que:

1. Ses fibres orbifoldes g\'en\'erales sont {\it sp\'eciales}.

2. Sa base orbifolde $(W\vert\Delta_W)$  est de type g\'en\'eral (ou un point si $(X\vert \Delta)$ est sp\'eciale).

Le coeur scinde donc canoniquement et fonctoriellement $(X\vert \Delta)$ en ses parties antith\'etiques sp\'eciale (les fibres orbifoldes) et de type g\'en\'eral (la base orbifolde).

Le ``coeur" domine (resp. est domin\'e par) toutes les fibrations sur $(X\vert \Delta)$ dont la base orbifolde est (resp. les fibres orbifoldes sont) de type g\'en\'eral (resp. sp\'eciales). Voir [Ca 04], [Ca 07] pour plus de d\'etails sur cette fibration.


\

  \label{crj} On montre aussi (conditionnellement en une version orbifolde $C^{orb}_{n,m}$ de la conjecture $C_{n,m}$ d'Iitaka) dans [Ca 07] que $c=(M\circ r)^n$ si $n=dim(X)$. Ce qui fournit une d\'ecomposition fonctorielle canonique du coeur en fibrations \'el\'ementaires $M$ et $r$, versions orbifoldes des fibrations de Moishezon-Iitaka et du ``quotient rationnel". Ceci montre que $X$ est sp\'eciale si (et seulement si, avec $C^{orb}_{n,m})$ c'est une tour de fibrations dont les fibres orbifoldes ont soit $\kappa=0$, soit $\kappa_+=-\infty$ (une version faible de la connexit\'e rationnelle). La condition ``$(X\vert \Delta)$ sp\'eciale" \'equivaut donc (conditionnellement) \`a la combinaison orbifolde de ces deux propri\'et\'es. Nous devons \`a nouveau renvoyer \`a [Ca 07] pour toutes ces notions.

\

 On va maintenant montrer que les ensembles de Green-Lazarsfeld et les quotients r\'esolubles du groupe fondamental d'une orbifolde K\"ahler coincident essentiellement (\`a un facteur ab\'elien fixe pr\`es) avec ceux de son ``coeur", ce qui permet, pour leur \'etude,  de se ramener au cas des orbifoldes de type g\'en\'eral, et en particulier, projectives. 



\begin{theorem}\label{redalg} Soit $c:(X\vert \Delta)\to (W\vert \Delta_W)$ le ``coeur"  de $(X\vert \Delta)$ , avec $(X\vert \Delta)$ lisse, et $X$ compacte K\" ahler. 

On note $G:=\pi_1(X\vert \Delta)$, et $G_W:=\pi_1(W\vert\Delta_W)$, avec $\Delta_W:=\Delta_{c,\Delta}$.

 Soit $\chi\in \Sigma_c^1(G)$, et $0\neq u\in Hom_G(DG_{ab},\bC_{\chi})$. Il existe alors $\chi'\in \Sigma_c^1(W)$ tel que: $\chi=c^*(\chi')$, et $v\in Hom_{G_W}((DG_W)_{ab},\bC_{\chi'})$ tel que: $u=c^*(v)$.

\end{theorem}

{\bf D\'emonstration:} L'assertion r\'esulte de \ref{glalb'}. En effet, puisque la fibration $h\circ g\circ f: (X\vert \Delta)\to (V\vert \Delta_V)$ (introduite au \S \ref{qralb}) a une base orbifolde qui est de type g\'en\'eral, il existe une factorisation de $h\circ g\circ f$ par $c$. On conclut alors par la fonctorialit\'e de $\Sigma^1_c$. L'existence de $v$ r\'esulte du lemme \ref{desc}.(3). $\square$

\

 L'\'etude des quotients r\'esolubles de  $G:=\pi_1(X\vert \Delta)$ est bas\'ee sur le lemme \ref{albalg} ci-dessous. Consid\'erons le diagramme commutatif: 

\

\
\centerline{
\xymatrix{(X\vert \Delta)\ar[r]^{f}\ar[d]_{c} &(Y\vert \Delta_Y)\ar[d]_{b}\\(W\vert \Delta_W)\ar[r]^d&(T\vert \Delta_T)\\
}}

\

\

dans lequel $f$ (resp. $d)$ est la factorisation de Stein du morphisme d'Albanese de $(X\vert \Delta)$ (resp. $(W\vert \Delta_W))$  introduit au \S\ref{albstein} ci-dessus, tandis que $c$ est le coeur de $(X\vert \Delta)$, et que $b$ est d\'eduite des morphismes d'Albanese de $X$ et $W$. Les structures orbifoldes sont les bases orbifoldes des fibrations du diagramme. Elles sont compatibles \`a la composition des fibrations, en choisissant des mod\`eles bim\'eromorphes ``nets" ad\'equats (voir [Ca 07, 3.14]).

\begin{lemma}\label{albalg} Dans le diagramme pr\'ec\'edent, on a les propri\'et\'es suivantes:



1. Le morphisme naturel: $c': X\to W\times_{T}Y$ est surjectif et connexe. 

2. $b:Y\to T$ est un fibr\'e principal de fibre $B:=Ker(Alb(X)\to Alb(V))$.

3. La fibre orbifolde g\'en\'erique $(Y_t\vert \Delta_{Y,t})$ de $b: (Y\vert \Delta_Y)\to (T\vert \Delta_T)$ est sans structure orbifolde (ie: $\Delta_{Y,t}:=\Delta_Y\cap Y_t=0)$.
\end{lemma}

{\bf D\'emonstration:} 1. La surjectivit\'e de $c'$ provient de ce que le morphisme d'Albanese des fibres g\'en\'eriques de $c$ est surjectif (d'apr\`es \ref{specalb}), puisque celles-ci sont sp\'eciales. La connexit\'e de $c'$ est \'evidente, puisque $c$ est connexe.  

2. L'application $b$ \'etant induite par un morphisme entre tores complexes, est un fibr\'e principal si ses fibres sont connexes (ce qui est \'evident, puisque celles de $b\circ f=d\circ c$ le sont), et sont des rev\^etements \'etales des translat\'es d'un sous-tore complexe de $Alb(X)$. Cette propri\'et\'e est vraie, puisque les fibres  lisses de $c$ sont sp\'eciales, et ne dominent donc aucun rev\^etement ramifi\'e d'un tore complexe, puisque leur morphisme d'Albanese est surjectif et connexe (par \ref{specalb}).

3. R\'esulte de l'assertion 2. de \ref{specalb}. $\square$

\

On notera (comme ci-dessus): $G_Y:=\pi_1(Y\vert \Delta_Y)$, $G_W:=\pi_1(W\vert \Delta_W)$, et $G_T:=\pi_1(T\vert \Delta_T)$.
Nous avons donc un diagramme comutatif de morphismes de groupes \`a lignes et colonnes exactes:

\
\centerline{
\xymatrix{
G\ar[r]^{f_*}\ar[d]_{c_*} &G_Y\ar[d]_{b_*}\\
G_W\ar[r]^{d_*}\ar[d]_{}&G_T{}\ar[d]_{}\\
1&1}}

\

\begin{corollary}\label{qralg} Avec ces notations, on a une suite exacte naturelle induite par $b$: $$1\to \pi_1(B)\to G_Y\to G_T\to 1$$ Cette suite est (non canoniquement) scind\'ee par des scindages induisant des produits de groupes: $G_Y\cong \pi_1(B)\times G_T$.
\end{corollary}

{\bf D\'emonstration:} On a (par [Ca 07, 11.2]) une suite exacte de groupes:$$1\to \pi_1(B\vert \Delta_B)\to G_Y\to G_T\to 1,$$ o\`u $(B\vert \Delta_B)$ est la restriction de $\Delta_Y$ \`a une fibre g\'en\'erique $F\cong B$ de $b$. Mais 
on a: $\Delta_B=0$ (par \ref{albalg} ci-dessus), et donc: $ \pi_1(B\vert \Delta_B)= \pi_1(B)$  

\

On a donc un diagramme commutatif \`a lignes et colonnes exactes:

\

\
\centerline{
\xymatrix{&1\ar[d]_{}&\\
\pi_1(B)\ar[r]^{}\ar[d]_{}&\pi_1(B)\ar[d]_{}&\\
G_Y\ar[r]^{}\ar[d]_{} & \pi_1(Y)\ar[r]^{}\ar[d]_{}&1\\
G_T\ar[r]^{}\ar[d]_{}&\pi_1(T)\ar[r]^{}\ar[d]_{}&1\\
1&1&}}

\
\

d'o\`u l'on d\'eduit ais\'ement l'assertion, en scindant la suite exacte $$1\to \pi_1(B)\to G_Y\to G_T\to 1$$ \`a l'aide de l'image inverse, not\'ee  $\widetilde{\pi_1(T)}$ dans $G_Y$ d'un groupe  $\pi_1(T)^* <\pi_1(Y)$, $\pi_1(T)^*$  isomorphe \`a $\pi_1(T)$ et induisant un isomorphisme: $\pi_1(Y)\cong \pi_1(T)^*\times \pi_1(B)$. 

Il est alors imm\'ediat de v\'erifier que $G_Y=\widetilde{\pi_1(T)}\ltimes \pi_1(B)$, que $\pi_1(B)$ est central dans $G_Y$, et que $\widetilde{\pi_1(T)}\cong G_T$. $\square$

\

Le r\'esultat suivant d\'ecoule alors imm\'ediatement de \ref{qralg} et de \ref{redalg}:

\begin{theorem}\label{qralg'} Soit $\varphi:G:=\pi_1(X\vert \Delta)\to R$ un quotient r\'esoluble de $G$. Soit (comme ci-dessus) $G_f:=Ker(f_*:G\to G_Y)$. Alors:

1. $\varphi(G_f)$ est un sous-groupe fini normal de $R$ (par \ref{qralb}). Soit $\rho:R\to R':=R/\varphi(G_f)$ le quotient, et $\varphi':=\rho\circ \varphi:G\to R'$.

2. Soit $A:=Ker(a_*:G\to G_W)$. Alors: $\varphi'(A)$ est un sous-groupe ab\'elien central de $R'$, isomorphe \`a un quotient de $\pi_1(B):=Ker(a_*:G_{ab}\to (G_W)_{ab})/Torsion$.

3. Pour tout entier $k\geq 0$, on a des \'equivalences de commensurabilit\'e: $$(G/D^{k+1}G)\equiv (G_Y/D^{k+1}G_Y)\cong (G_T/D^{k+1}G_T)\times \pi_1(B)\equiv (G_W/D^{k+1}G_W)\times \pi_1(B),$$ o\`u: $H\equiv L$ signifie que les groupes $H$ et $L$ sont isomorphes, apr\`es passage a quotient par des sous-groupes normaux finis ad\'equats.

Autrement dit: les quotients r\'esolubles canoniques de $G$ sont (\`a un groupe fini fixe pr\`es) des produits de ceux de $G_W$ (groupe fondamental de la 
de la base orbifolde du ``coeur" de $(X\vert \Delta))$ par des quotients de $\pi_1(B)$. \end{theorem}

{\bf D\'emonstration:} Les assertions 1 et 2 sont imm\'ediates, compte-tenu de \ref{rqralb} et de \ref{qralg}. Assertion 3: la premi\`ere \'equivalence r\'esulte de \ref{rqralb} appliqu\'ee au quotient naturel $\varphi: G\to R:=G/D^{k+1}G$. Remarquer en effet que $\varphi(G_f)$ est fini et normal. On a donc une factorisation $\varphi':G_Y/D^{k+1}G_Y\to R/\varphi (G_f):=R'$, qui montre la premi\`ere \'equivalence annonc\'ee. La seconde d\'ecoule de \ref{qralg}. La troisi\`eme n'est autre que la premi\`ere, appliqu\'ee \`a $(W\vert \Delta_W)$, apr\`es produit avec $\pi_1(B)$. $\square$

\begin{remark}\label{rkalbalg} Une r\'eponse positive \`a la question \ref{qqralb}.1 pour le groupe $\varphi(G_{gf})$ permettrait, \`a l'aide des arguments de remplacer $Y$ et $T$ par les images d'Albanese de $X$ et $W$ dans les \'enonc\'es pr\'ec\'edents. Pour les quotients nilpotents, ces \'enonc\'es concernant les images d'Albanese de $X$ et $W$ sont vrais, par [Ca 95]. \end{remark}

\section{Quotients nilpotents des groupes de K\" ahler.}

\subsection{Rappels sur quotients et groupes nilpotents.}\label{rgn}

Soit $G$ un groupe, et $s\geq 1$ un entier. On note $C^s(G)$ le $s$-i\`eme terme de sa suite centrale descendante, d\'efinie par r\'ecurrence par: $C^1(G):=G,$ et $ C^{s+1}(G):=[G,C^s(G)]$, le groupe engendr\'e par les commutateurs de $G$ et de $C^s(G)$ si $s>1$. Pour tout $s\geq 1$, $C^s(G)$ est un sous-groupe caract\'eristique de $G$.

Le groupe $G$ est dit {\it nilpotent de classe au plus s} si $C^{s+1}(G)$ est trivial (ie: r\'eduit \`a son \'el\'ement neutre). On note alors $\nu (G)$ le plus petit entier $s\geq 1$ tel que $C^{s+1}(G)=\{1\}$. On appelle $\nu(G)$ la classe de nilpotence ({\it exacte}) de $G$.
Donc $\nu(G)=0$ (resp $1$) si et seulement si $G$ est trivial (resp. ab\'elien non trivial).

Une extension $1\to Z\to G\to H\to 1$ d'un groupe nilpotent $H$ par un groupe ab\'elien $Z$ est nilpotente si elle est {\it centrale} (ie: si $Z\subset Z(G)$, le centre de $G$).

On notera $q_s: G\to G/C^{s+1}(G):=G_s$ le quotient naturel: $G_s$ est ainsi le plus grand quotient nilpotent de $G$ de classe au plus $s$ (les autres \'etant des quotients de $G_s$). En g\'en\'eral, $\nu(G_s)<s$ (par exemple, si $G$ est ab\'elien, alors $\nu(G_s)=1$ pour tout $s>1$).

La plupart des propri\'et\'es de finitude des groupes ab\'eliens s'\'etendent aux  groupes nilpotents (mais non pas aux groupes r\'esolubles). Si $G$ est un groupe nilpotent, alors (voir [Ro] pour les 4 premi\`eres propri\'et\'es, et [R] pour la derni\`ere):

1. L'ensemble des \'el\'ements de torsion $T(G)$ est un sous-groupe (caract\'eristique) de $G$. On note $\tau_G:G\to G^*:=G/T(G)$ le quotient. Alors $G^*$ est sans torsion. On notera $\nu^*(G):=\nu(G^*)\leq \nu(G)$.

2. Si $G$ est de type fini, $T(G)$ est fini.

3. Si $G$ est de type fini, il est de pr\'esentation finie.

4. Si $G$ est de type fini, il est lin\'eaire (ie: plongeable dans un $Gl(n,\bR)$). Plus pr\'ecis\'ement:

5. Si $G$ est nilpotent, sans torsion, et de type fini, il existe un unique groupe de Lie r\'eel et simplement connexe $G_{\bR}$, d\'efini sur $\bQ$, tel que $G$ soit un r\'eseau cocompact de $G_{\bR}$. (On peut d\'efinir $G_{\bQ}$ intrins\`equement comme la r\'eunion des ensembles des racines $m-$i\`emes (uniques, car $G$ est nilpotent sans torsion) des \'el\'ements de $G$, pour $m>0$ entier). Le groupe $G_{\bR}$ est la {\it compl\'etion de Mal\c cev} de $G$. On notera alors $LG$ son alg\`ebre de Lie (d\'efinie sur $\bQ$). On a donc aussi, pour tout $s>0$: $L(C^s(G))=C^sL(G)$, en notant $C^sL(G)$ le $s$-i\`eme terme de la suite descendante d'une alg\`ebre de Lie $L(G)$.

Si $G$ est un groupe arbitraire de type fini, alors $G^*_s:=(G_s)^*$ est aussi de type fini, nilpotent sans torsion et de classe de nilpotence au plus $s$. On notera donc $LC^s(G^*_s):=LC^s(G)$. Remarquer que $LC^s(G^*_{m+s})=LC^s(G^*_s),$ pour tout entier $m\geq 0$. 

D'apr\`es la propri\'et\'e 1 ci-dessus, si $G$ est un groupe arbitraire, et $s>0$ entier, on peut donc d\'efinir la compos\'ee $q^*_s:=(\tau_{G_s})\circ q_s:G\to G^*_s$, et $C_*^{s+1}(G):=Ker(q^*_s)$, que nous appellerons (abusivement) $(s+1)$-i\`eme terme de la suite entrale descandante {\it sans torsion} de $G$.

Les constructions pr\'ec\'edentes sont \'evidemment fonctorielles en $G$. Si $f:G\to H$ est un morphisme de groupes, il induit donc des morphismes $f_s:C^s(G)\to C^s(H)$, $f^s:G_s\to H_s$, $f^*_s:G^*_s\to H^*_s$, surjectifs si $f$ l'est. Et aussi, si $G,H$ sont de type fini et sans torsion: $f_{\bR}:G_{\bR}\to H_{\bR}$, ainsi que $Lf:LG\to LH$, qui peuvent \^etre compos\'es avec les pr\'ec\'edents.

\subsection{Morphismes stricts de groupes}

\

\begin{definition}\label{defstr} Soit $f:G\to H$ un morphisme de groupes, $G$ et $H$ \'etant de types finis. On dit que $f$ est un {\bf morphisme strict} si, pour tout entier $j\geq 1$, le morphisme induit $Lf_s^*:LC^s(G)\to LC^s(H)$ est strict, c'est-\`a-dire, si: $f(LC^s(G))\cap LC^j(H)=f(LC^j(G))$, pour tout $0\leq s\leq j$. (En g\'en\'eral, le second membre est seulement inclus dans le premier).
\end{definition}

Si $G=K\times H$ est le produit direct de deux groupes, et si $f$ est la seconde projection, alors $f$ est un morphisme strict. Cette notion ne semble pas avoir de caract\'erisation ou de propri\'et\'es de stabilit\'e de formulation simple en termes alg\'ebriques. Par contre, la g\'eom\'etrie K\" ahl\'erienne en fournit de nombreux exemples (qui inspirent la terminologie):

\begin{theorem}\label{str} Soit $g:X\to Y$ une application holomorphe entre vari\'et\'es K\" ahl\'eriennes compactes et connexes, et $f=g_*:\pi_1(X)\to \pi_1(Y)$ le morphisme induit (on omet les points-base). Alors $f$ est un morphisme de groupes strict.
\end{theorem}

{\bf D\'emonstration:} C'est une cons\'equence directe de deux r\'esultats profonds. Le premier, d\^u \`a R. Hain ( voir [H 87], et [P-S, chap. 8]), affirme que, pour tout entier $s>0$, l'alg\`ebre de Lie $LC^s(\pi_1(X))$ est munie d'une structure de Hodge mixte naturelle, fonctorielle en $X$, avec pour filtration par le poids: $W^j(LC^s(\pi_1(X)))=LC^j(\pi_1(X))$, pour tout $j\geq 0$. Le second r\'esultat, d\^u \`a P. Deligne ([D 71], voir aussi [P-S, corollary 3.6, p. 65]) affirme que tout morphisme $f: L\to M$ de structures de Hodge mixtes est strict, et en particulier, que $f(L)\cap W^j(M)=f(W^j(L))$, pour tout entier $j\geq 0$. $\square$

\

Nous utiliserons la notion de morphisme strict de groupes dans la situation suivante.

 Soit $h:K\to G$, et $f:G\to Q$ des morphismes de groupes de type fini tels que la compos\'ee $K\to G\to Q\to \{1\}$ soit exacte. Soit alors $s>0$ un entier, et $q^*_s:G\to G^*_s$ le quotient naturel. 
 
 \
 
 On en d\'eduit comme suit un diagramme commutatif de groupes, \`a lignes exactes:
 
 \

 \
\centerline{
\xymatrix{ K\ar[r]^{}\ar[d]_{}&G\ar[r]^{}\ar[d]_{}&Q\ar[r]^{}\ar[d]_{}&1\\K'\ar[r]^{}& G^*_s\ar[r]^{}&Q'\ar[r]^{}&1\\
}}

\

Soit $K':=K/K\cap Ker (q^*_s)$, et $h_s:K\to K'$ les quotients d\'eduits de $h$ et $K$. Le groupe $K'$ est naturellement un sous-groupe distingu\'e (ou normal) de $G^*_s$. Soit $G^*_s\to G^*_s/K':=Q'$ le quotient. Il existe alors un unique morphisme de groupes (surjectif) $Q\to Q'$ rendant commutatif le diagramme pr\'ec\'edent.

Les groupes $K'$ et $Q'$ sont donc nilpotents de classe au plus $s$, et $K'$ est sans torsion.

\begin{proposition}\label{str'} Dans la situation pr\'ec\'edente, on suppose que $h:K\to G$ est strict. Alors:

 $$\nu(G^*_s)=max\{\nu(K'),\nu^*(Q')\}:=m$$
\end{proposition}

{\bf D\'emonstration:} On a, bien s\^ur, toujours: $\nu(G^*_s)\geq m$ sans hypoth\`ese sur $h$. On peut supposer que $Q'$ est sans torsion (en effet: remplacer $Q'$ par $Q":=Q'/T(Q')$; le noyau du morphisme compos\'e $G^*_s\to Q"$ est alors un sous-groupe $K"$ (nilpotent sans torsion) de $G^*_s$ dans lequel $K'$ est d'indice fini. On a donc: $\nu(K")=\nu(K')$. Il suffit alors de consid\'erer la suite exacte $K"\to G^*_s\to Q"\to \{1\}$, puisque les alg\`ebres de Lie des compl\'etions de Mal\c cev de $K'$ et $K"$ coincident, de sorte que le morphisme $K"\to G^*_s$ est \'egalement strict).

Pla\c cons-nous dans $L^s(G)=L^s(G^*_s)$. 

(1) On a: $L^{m+1}(K")=L(K")\cap L^{m+1}(G^*_s)=\{0\}$, puisque $h$ est strict et $m\geq \nu(K")$.

(2) D'autre part, $L^sh(L^{m+1}(G^*_s))\subset L^{m+1}(Q")=\{0\}$, puisque $m\geq \nu(Q")$.

(3) Enfin, $L(K")=Ker(L^sh)$, puisque $K"=Ker(G^*_s\to Q")$. 

On d\'eduit de (2) et (3) que $L^{m+1}(G^*_s)\subset L(K")$, puis de (1) que $L^{m+1}(G^*_s)=\{0\}$. 

Puisque $L^{m+1}(G^*_s)=LC^{m+1}(G^*_s)$, on a donc bien: $C^{m+1}(G^*_s)=\{1\}$ $\square$

\

\subsection{Fibrations et classe de nilpotence des groupes de K\" ahler.}

\

 Soit $X,Y$ des vari\'et\'es K\"ahl\'eriennes compactes et connexes, et $g:X\to Y$ une {\it fibration}, c'est-\`a-dire: une application holomorphe surjective et \`a fibres connexes. Soit $X_y$ une fibre lisse arbitraire de $f$. L'inclusion de $X_y$ dans $X$ induit un morphisme de groupes $h:K:=\pi_1(X_y)\to G:=\pi_1(X)$, dont l'image est un sous-groupe normal de $\pi_1(X)$. Le morphisme de groupes $g_*=f:G:=\pi_1(X)\to Q:=\pi_1(Y)$ est surjectif.

 La compos\'ee $g_*\circ h:\pi_1(X_y)\to \pi_1(Y)$ est triviale, mais la suite de groupes $K:=\pi_1(X_y)\to G:= \pi_1(X)\to Q:=\pi_1(Y)\to \{1\}$ n'est pas exacte en $\pi_1(X)$ en g\'en\'eral (en pr\'esence de fibres multiples).
 
 Consid\'erons (comme dans la proposition \ref{str'}) le quotient: $q^*_s:G\to G^*_s:=\pi_1(X)^*_s$, puis le noyau $\sK$ du morphisme compos\'e $q^*_s\circ h:K\to G^*_s$, ensuite le quotient $K':=K/\sK$, et enfin le quotient $G^*_s\to Q':=G^*_s/K'$ (qui n'est pas {\it a priori} un quotient de $Q$, mais seulement de $G/K$). 

\begin{theorem}\label{fibr} Dans la situation pr\'ec\'edente, les trois propri\'et\'es suivantes sont satisfaites:

1. $Q'$ est un quotient de $Q$.

2. $\nu(G^*_s)=max\{\nu(K'),\nu^*(Q')\}$. 

3. $\nu(G^*_s)\leq max\{\nu(K^*_s),\nu(Q^*_s)\}$.
\end{theorem}

{\bf D\'emonstration:} La propri\'et\'e (3) d\'ecoule de (2) et (1), puisque $K'$ (resp. $Q'$) est un quotient de $K^*_s$ (resp. $Q_s$). La propri\'et\'e (2) r\'esulte de \ref{str} et de \ref{str'}. Enfin, le lemme \ref{q} suivant implique aussi la propri\'et\'e (1), puisque $Q'$ est nilpotent de classe au plus $s$ $\square$

\begin{lemma}\label{q} Avec les notations $K=\pi_1(X_y)$, $G=\pi_1(X)$ pr\'ec\'edentes, si $\bar Q:=G/K$, alors pour tout entier $s>0$, le quotient naturel $\bar Q\to Q$ induit des isomorphismes $\bar Q^*_s\cong Q^*_s$. Plus g\'en\'eralement, tout morphisme de groupes $\bar Q\to H$ dont l'image est sans torsion se factorise par $Q$.
\end{lemma}

{\bf D\'emonstration:}Il suffit d'\'etablir la seconde assertion, puisque les groupes $\bar Q^*_s$ sont lin\'eaires sans torsion. Nous allons utiliser pour cela la base orbifolde $\Delta_g$ de la fibration $g$ (voir [Ca 07]]). Les assertions de l'\'enonc\'e sont invariantes par transformation bim\'eromophe de $X$. Nous pouvons donc supposer la fibration $g: X\to Y$ {\it nette} (voir [Ca 07,chap. 3.2]) et telle que sa base orbifolde $(Y\vert\Delta_g)$ soit lisse (ie: $Y$ est lisse, et le support du diviseur $\Delta_g$ est \`a croisements normaux). 

Dans ce cas, $\bar Q=\pi_1(Y\vert \Delta_g)$ ([Ca 07,chap. 11.2], par exemple), et le noyau du quotient $\bar Q:=\pi_1(Y\vert \Delta_g)\to Q=\pi_1(Y)$ est engendr\'e par des \'el\'ements de torsion. Le lemme en r\'esulte, puisque ces \'el\'ements ont une image triviale dans $H$ $\square$


\begin{corollary}\label{fibr'} Soit $g:X\to Y$ une fibration, avec $X$ K\"ahler compacte et connexe\footnote{On doit supposer $g$ ``nette" (voir [Ca 07]). Cette condition est r\'ealis\'ee apr\`es modifications de $X$ et $Y$.}. Supposons que $K:=\pi_1(X_y)_X:=Im(\pi_1(X_y)\to \pi_1(X))$ et $Q:=\pi_1(Y\vert \Delta)$ sont r\'esolubles. Alors $\pi_1(X), K$ et  $Q$ sont virtuellement nilpotents, et $\nu(G^*)=max\{\nu(K^*),\nu(Q^*)\}$.

En particulier, si $K$ et $Q$ sont presque-ab\'eliens, $G$ l'est aussi.\end{corollary}

\begin{example} On d\'eduit ais\'ement de \ref{fibr'} que si $f:X\to Y$ est une submersion propre et connexe sur un tore complexe, et si les fibres de $f$ sont des tores complexes, alors $X$ est rev\^etue par un tore complexe si $X$ est K\"ahler (propri\'et\'e sugg\'er\'ee par une observation similaire du referee). Remarquer en effet que la fibration est alors isotriviale. \end{example}

\begin{re} 

\

1. Le th\'eor\`eme \ref{fibr} s'\'etend ais\'ement au cadre des orbifoldes \`a multiplicit\'es finies (gr\^ace au lemme \ref{q} ci-dessus). Par contre le cas logarithmique n\'ecessiterait des ingr\'edients nouveaux (tels que ceux de [M 78]).

2. Il serait int\'eressant de savoir si la propri\'et\'e (2) du th\'eor\`eme \ref{fibr} subsiste si l'on y remplace $G^*_s=\pi_1(X)^*_s$ par un quotient nilpotent sans torsion arbitraire de $G=\pi_1(X)$. 

3. Observons que cette propri\'et\'e ne subsiste pas si l'on remplace $K$ par un sous-groupe normal ne provenant pas d'une fibration. On le voit sur l'exemple de groupe de K\" ahler  Heisenberg $G$ construit dans [Ca 95]: on a en effet une suite exacte: $$1\to \Bbb Z\to G\to A\oplus B\to 1,$$ avec: $A\cong B\cong \Bbb Z^{\oplus 2m}, m\geq 4$. L'image r\'eciproque $A'$ de  $A$  dans $G$ est ab\'elien, ainsi que $G/A'\cong B$, alors que $2=\nu(G)>max\{\nu(A'),\nu(G/A')\}=1$. \end{re}

\subsection{Question.}\label{qq}

Nous venons de voir que les quotients des sous-groupes d'indice fini de $G=\pi_1(X)$, $X$ copacte lisse et  K\"ahler sont des groupes de type fini arbitraires  si et seulement si $X$ admet (apr\`es rev\^etement fini \'etale $X'$ ad\'equat) un morphisme sur une courbe hyperbolique. Il est donc essentiel de pouvoir caract\'eriser cette situation aussi simplement que possible. Nous avons vu qu'elle \'equivaut aux conditions: $\Sigma^1(G')$ infini, ou encore $DG'$ de type infini.

On propose ci-dessous deux autres conditions potentiellement \'equivalentes:

\begin{question} Soit $X$ compacte et K\" ahler connexe, et $X'$ un rev\^etement \'etale fini ad\'equat, avec $G':=\pi_1(X')$. A-t'on \'equivalence entre les propri\'et\'es suivantes:

1. Il existe $X'$ admettant une fibration $f':X'\to C'$ sur une courbe de genre $g\geq 2$.

2. Pour tout $N>0$, il existe $X'$ tel que $q(X'):=h^{1,0}(X')\geq N$.

3. Il existe $X'$ tel que $\nu((G')^*_s)=s, \forall s>0$, avec les notations de \ref{rgn}.
\end{question}

La condition 1 implique les deux autres. L'implication $2\Longrightarrow 1$ pourrait \^etre \'eventuellement \'etablie en montrant que $DG'$ n'est alors pas de type fini.  L'implication $3\Longrightarrow 2$ semble n\'ecessiter l'utilisation des ensembles de Green-Lazarsfeld et la structure de Hodge mixte de Hain sur les quotients nilpotents de $G'$.

\section{Bibliographie}

[A 92] D. Arapura. Higgs line bundles, Green-Lazarsfed sets, and maps of K\" ahler manifolds to curves. Bull. AMS. 26 (1992), 310-314.

[AN 99] D. Arapura-M. Nori. Solvable fundamental groups of algebraic varieties and K\" ahler manifolds. Comp. Math. 116 (1999), 173-188.

[Be 92] A. Beauville. Annulation du $H^1$ pour les fibr\'es en droites plats. LNM 1507 (1992), 1-15.

[B-S] R. Bieri-R. Strebel. Geometric invariants for discrete groups. Cours polycopi\'e. L'article original est [B-N-S 87].

[B-N-S 87] R. Bieri-W. Neumann-R. Strebel. A geometric invariant of discrete groups. Inv. Math. 90 (1987), 451-477.

[Bru 99] A. Brudnyi. Classification theorem for a class of flat connections and representations of K\" ahler groups. Michigan J. Math. 46 (1999), 489-514.

[Ca 85] F. Campana. R\'eduction d'Albanese d'un morphisme K\" ahl\'erien propre. I, II. Compositio Mathematica 54 (1985), 373-416.

[Ca 95] F. Campana. Remarques sur les groupes de K\" ahler nilpotents. Ann. Sc. ENS. 28 (1995), 307-316.

[Ca 98] F. Campana. Connexit\'e ab\'elienne des vari\'et\'es K\" ahl\'eriennes compactes. Bull. SMF 126 (1998), 483-506.

[Ca 01] F. Campana. Ensembles de Green-Lazarsfeld et quotients r\'esolubles des groupes de K\" ahler. J. Alg. Geom. 10 (2001), 599-622.

[Ca 04] F. Campana. Orbifolds, special varieties and classification theory. Ann. Inst. Fourier 54 (2004), 499-635.

[Ca 07] F. Campana. Orbifoldes sp\'eciales et classification bim\'eromorphe des vari\'et\'es K\" ahl\'eriennes compactes. arXiv 0705.0737.

[Cat 91] F. Catanese. Moduli and classification of irregular K\" ahler Manifolds with Albanese general type fibrations. Inv. Math. 104 (1991), 263-289 . 

[Cat 94] F. Catanese. Fundamental groups with few relations. Proceedings of the School-Conference ``Higher Dimensional Complex Varieties. Trento, June 1994. Walter De Gruyter, Berlin (1996), 163-165.

[Cat 07] F. Catanese. Differentiable and Deformation type of Algebraic Surfaces, Real and Symplectic Structures. in Symplectic 4-Manifolds and Algebraic Surfaces. Proceedings of the Cetraro Conference, September 2007. SpringerLecture Notes (to appear), 55-168.

[Cl 08] B. Claudon. Invariance de la Gamma-dimension pour certaines familles K\" ahl\'eriennes de dimension 3. A para\^itre au Math. Z. Arxiv: 0802.2894.

[Cl 08']B. Claudon.$\Gamma$-reduction for smooth orbifolds. arXiv:0802.2894.
 
[D 71]P. Deligne. Th\'eorie de Hodge II. Publ. IHES 40 (1971), 5-58.

[D-M 93] P. Deligne-D. Mostow. Commensurabilities among lattices in PU(1,n). Annals of Mathematic Studies 132. Princeton University Press 1993.

[De06] T. Delzant. L'invariant de Bieri-Neumann-Strebel des groupes fondamentaux des vari\'et\'es K\" ahl\'eriennes. arXiv. math/0603038. (A para\^itre aux Math. Ann.)

[De07] T. Delzant. Trees, valuations, and the Green-Lazarsfeld set. arXiv. math/0702477. A para\^itre \`a GAFA.

[DG 05] T. Delzant, M. Gromov. Cuts in K\" ahler groups. Progress in Mathematics 248. Birkh\" auser 2005.

[F 83] A. Fujiki. On the structure of a complex manifold in C. Adv. Studies in Pure Math. 1 (1983), 231-302.

[GL 87] M. Green-R. Lazarsfeld. Deformation theory, generic vanishing theorems, and some conjectures of Enriques, Catanese and Beauville. Inv. Math. 90 (1987), 389-407.

[H 87] R. Hain. The De Rham homotopy of complex algebraic varieties I and II. K-theory 1 (1987), 271-324 and  481-497.

[K 81] O. Kharlampovich. A finitely presented solvable group with unsolvable word problem. Izv. Akad. Nauk. ser. Mat. 45 (1981), 852-873.

[Ko 93] J. Koll\`ar. Shafarevich maps and plurigenera of algebraic varieties. Inv. Math.113 (1993), 177-215.

[M 78] J. Morgan. The algebraic topology of smooth algebraic varieties. Publ. Math. IHES 48 (1978), 137-204. Correction: Publ. Math. IHES 64 (1986), p. 185

[N 87] M. Namba. Branched coverings and algebraic functions. Pitman Research Notes in Mathematics Series, vol. 161. Longman Scientific and Technical (1987).

[P-S 08] C. Peters-J. Steenbrink. Mixed Hodge Structures. Erg. der Math. und ihrer Grenzgebiete. 3. Folge. Vol. 52. Springer Verlag  (2008).

[R] M.Ragunathan. Discrete subgroups of Lie groups. Ergebnisse der Mathematik und Ihrer Grenzgebiete. Springer Verlag. 

[Ro]D. Robinson. A course in the theory of groups. GTM 80. Springer Verlag (1996).

[Se 68]J.P. Serre. Corps locaux. Hermann 1968.

[Si 93] C. Simpson. Subspaces of moduli spaces of rank one local systems. Ann. Sci. ENS. 26 (1993), 361-401.

[Si 93-2] C. Simpson. Lefschetz theorem for the integral leaves of a holomorphic one-form. Compositio Math. 87 (1993), 99-113.

[U 75] K. Ueno. Classification theory of algebraic varieties and compact complex spaces. LNM 439. Springer Verlag 1975.

\

{\bf Adresse}

Fr\'ed\'eric Campana

Universit\'e Nancy 1.

D\'epartement de Math\'ematiques. 

e-mail: campana@iecn.u-nancy.fr

\end{document}